\documentclass[11pt,a4paper]{article}
\usepackage{graphicx}
\usepackage{amsmath}
\usepackage{geometry}
\usepackage{amsfonts}
\usepackage{amssymb}
\newtheorem{theorem}{Theorem}[section]

\newtheorem{corollary}{Corollary}[section]

\newtheorem{lemma}{Lemma}[section]

\newtheorem{proposition}{Proposition}[section]
\newtheorem{remark}{Remark}[section]

\newenvironment{proof}[1][Proof]{\textbf{#1.} }{\ \rule{0.5em}{0.5em}}
\numberwithin{equation}{section}
\begin{document}

\title{Critical behavior of the massless free field at the depinning transition}
\author{Erwin~Bolthausen
\and Yvan~Velenik}
\date{\today}
\maketitle
\begin{abstract}
We consider the $d$-dimensional massless free field localized by a $\delta
$-pinning of strength $\varepsilon$. We study the asymptotics of the variance
of the field (when $d=2$), and of the decay-rate of its 2-point function (when
$d\geq2$), as $\varepsilon$ goes to zero, for general Gaussian interactions.
Physically speaking, we thus rigorously obtain the critical behavior of the
transverse and longitudinal correlation lengths of the corresponding
$d+1$-dimensional effective interface model in a non-mean-field regime. We also
describe the set of pinned sites at small $\varepsilon$, for a broad class of
$d$-dimensional massless models.
\end{abstract}

\section{Introduction}

The behavior of a two-dimensional interface at phase transitions has been much
studied in the physics literature, especially regarding some models of
wetting. The latter problem arises when one considers an interface above an
attractive wall. Then there is a competition between attraction by the wall
and repulsion due to the decrease of entropy for interfaces close to it.
Often, tuning some external parameter (the temperature, or the strength of the
attraction), two behaviors are possible: either energy wins, and the interface
stays localized along the wall, or entropy wins, and the interface is repelled
at a distance from the wall diverging as the size of the system grows. The
corresponding transition is called wetting transition. Usually in Nature this
transition is first-order, which means here that the average height of the
interface above the wall stays uniformly bounded as the parameter approaches
the critical value from the localized phase, and makes a jump ``to infinity''
(in the thermodynamic limit) at the transition. There are however cases when
this transition is second-order (the two-dimensional Ising model is a nice
theoretical example, but this behavior can also be observed in real systems);
this is the so-called critical wetting.
In this case, the average height of the interface diverges continuously as the
critical value is approached. It is then of interest to characterize this
divergence. We refer to~\cite{FJ} for references to the (non-rigorous) results
which have been obtained.

Unfortunately, very little is known rigorously about the behavior of
two-dimensional interfaces at a critical wetting transition, even for simple
effective interface models. There are some results on part of the so-called
``mean-field'' regime~\cite{DMR, Lemberger}, but nothing concerning the more
interesting ones.

In the present work, we study the critical behavior of a $d$-dimensional
interface localized by a $\delta$-pinning (defined below). The main focus will
be on the most difficult and physically most relevant two-dimensional case,
but the other cases will also be discussed. Though this problem is
clearly simpler than the wetting transition, it has the advantage of being
non-mean-field, while being rigorously tractable; we make some additional
comments on the wetting problem at the end of Section~\ref{SectResults}.

Let $\Lambda\Subset\mathbb{Z}^{d}$. We consider the following class of
massless gradient models in $\Lambda$, with $0$-boundary conditions described
by the following probability measures on $\mathbb{R}^{\mathbb{Z}^{d}}$
($\delta_{0}$ is the point mass at $0$):
\begin{equation}
\mu_{\Lambda}\left(  \mathrm{d}\phi\right)  \overset{\mathrm{def}}{=}\frac
{1}{Z_{\Lambda}}\exp\left[  -\frac{\beta}{2}\sum_{x,y}p\left(  x-y\right)
V\left(  \phi_{x}-\phi_{y}\right)  \right]  \prod_{x\in\Lambda}\mathrm{d}%
\phi_{x}\,\prod_{x\not \in\Lambda}\delta_{0}\left(  \mathrm{d}\phi_{x}\right)
,\label{free}%
\end{equation}
where $V$ is an even and convex function, and $\beta>0.$ We assume that
$p(x)=p(-x)\geq0,$ $\sum_{x\in\mathbb{Z}^{d}}p(x)=1,$ for any $x\in
\mathbb{Z}^{d}$ there exists a path $0\equiv x_{0},x_{1},\ldots,x_{n}\equiv x$
such that $p(x_{k}-x_{k-1})>0,\,k=1,\ldots,n,$ and at least
\begin{equation}
\sum_{x\in\mathbb{Z}^{2}}p(x)\left|  x\right|  ^{2+\delta}<\infty
\label{CondH1}%
\end{equation}
for some $\delta>0.$

We denote by $\mu_{\Lambda}^{\star}$ the (Gaussian) measure corresponding to
the particular choice $V(x)=\tfrac{1}{2}x^{2}$. A $\star$ superscript will
always be used for quadratic interactions.

It is well-known that for $d=2$ these measures describe a random field with
unbounded fluctuations as $\Lambda\nearrow\mathbb{Z}^{d},$ diverging
logarithmically with the size of $\Lambda$ if the limit is taken along a
sequence of cubes, say, while for $d\geq3$ the variance stays bounded. For the
Gaussian case this follows from the well known random walk representation of
the covariances
\begin{equation}
\mu_{\Lambda}^{\ast}\left(  \phi_{x}\phi_{y}\right)  =\frac{1}{\beta
}\mathrm{E}_{x}\left(  \sum_{n=0}^{\tau_{\Lambda}}\mathrm{I}\left(
X_{n}=y\right)  \right)  ,\label{RW}%
\end{equation}
where $\left(  X_{n}\right)  $ is a random walk, starting at $x$ under
\textrm{P}$_{x}$, with transition probabilities given by
\textrm{P}$_{x}\left(
X_{1}=y\right)~=~p\left(  y-x\right)$, $\tau_{\Lambda}$ is the first exit
time from $\Lambda,$ and $\mathrm{I}\left(  \cdot\right)  $ denotes the
indicator function of a set. A two-dimensional symmetric random walk
satisfying (\ref{CondH1}) is recurrent, and so for $d=2$ the divergence of the
variances as $\Lambda\nearrow\mathbb{Z}^{2}$ follows. In higher dimensions,
random walks are transient, and therefore, the variance stays bounded. Notice
however, that even in two dimensions, a random walk satisfying $\sum
_{x}p(x)\left|  x\right|  ^{2-\delta}=\infty$ for some $\delta>0$ is transient.

For more general convex interaction functions $V,$ the corresponding results
follow by an application of the Brascamp-Lieb inequality (see \cite{BrLi}).

It turns out, however, that the addition of an arbitrarily weak self-potential
breaking the continuous symmetry of the Hamiltonian, $\phi\rightarrow\phi+c$,
$c\in\mathbb{R}$, is enough to localize the field. More precisely, if $a$ and
$b$ are two strictly positive real numbers, then we perturb the measures by
modifying them with a ``square well'' potential:
\begin{equation}
\mu_{\Lambda}^{a,b}(\,\cdot\,)\overset{\mathrm{def}}{=}\frac{\mu_{\Lambda
}\,\left(  \,\cdot\;\exp[b\sum_{x\in\Lambda}\mathrm{I}\left(  \left|  \phi
_{x}\right|  \leq a\right)  ]\right)  }{\mu_{\Lambda}\left(  \exp[b\sum
_{x\in\Lambda}\mathrm{I}\left(  \left|  \phi_{x}\right|  \leq a\right)
]\right)  }.\label{DefSquarewell}%
\end{equation}
Another type of pinning, mathematically slightly more convenient, has also
been investigated, the so-called $\delta$-pinning. It corresponds to the weak
limit of the above measures when $a\rightarrow0$ and $2a(e^{b}-1)=\varepsilon
$, for some $\varepsilon>0$, and has the following representation:
\begin{equation}
\mu_{\Lambda}^{\varepsilon}(\mathrm{d}\phi)=\frac{1}{Z_{\Lambda}^{\varepsilon
}}\,\exp\left[  -\tfrac{\beta}{2}\sum_{x,y}p(x-y)\,V(\phi_{x}-\phi
_{y})\right]  \,\prod_{x\in\Lambda}\left(  \mathrm{d}\phi_{x}+\varepsilon
\delta_{0}(\mathrm{d}\phi_{x})\right)  \prod_{x\not \in\Lambda}\delta
_{0}(\mathrm{d}\phi_{x}).\label{DefDeltapin}%
\end{equation}
The most natural question in two dimensions is if a thermodynamic limit as
$\Lambda\nearrow\mathbb{Z}^{2}$ of these measures exists. The answer is most
probably ``yes'', but we cannot prove this, except in the Gaussian case with
$\delta-$pinning (Proposition \ref{prop_thermolimit} below). A somewhat
simpler question is whether the variance stays bounded uniformly in $\Lambda.$ This
was shown for $\mu_{\Lambda}^{a,b}$ in the Gaussian nearest neighbor case in
\cite{DMRR}, and was finally proved in \cite{DV} much more generally, assuming
only $V^{\prime\prime}\geq\mathrm{const.}>0. $ Moreover, it was shown in
\cite{IV} that the covariances $\mu_{\Lambda}^{\varepsilon}\left(  \phi
_{x}\phi_{y}\right)  $ decay exponentially in $\left|  x-y\right|  ,$
uniformly in $\Lambda,$ provided $0<\mathrm{const.}\leq V^{\prime\prime}%
\leq\mathrm{const.}<\infty$ (see also \cite{BB} for the Gaussian nearest
neighbor case). The discussion in \cite{IV} is restricted to the $\delta
$-pinning case, but it could probably be extended to the square well case at
least for quadratic interactions.

The aim of the present paper is to obtain a precise description of the
behavior of the variance of the field (or equivalently, in a more physical
terminology, of the transverse correlation length) and of the rate of decay of
the covariance (or of the longitudinal correlation length), as one approaches the
depinning transition, i.e. as the strength $\varepsilon$ of the pinning
potential goes to zero. The latter question is also of interest for $d\geq3.$
For the Gaussian $\delta$-pinning case, we determine exactly the divergence of
the variance for $d=2$ (Theorem \ref{thm_critical}) as a function of the
pinning parameter $\varepsilon,$ and the $\varepsilon$-dependence of the mass
for $d\geq2$, including the correct power of the logarithmic correction for
$d=2$ to the power law dependence in $\varepsilon$ (Theorem
\ref{thm_critical2}).

There are two main ingredients to our approach. By a simple expansion like
expanding the product $\prod_{x\in\Lambda}$ in (\ref{DefDeltapin}), we obtain
a representation of the random field as a mixture of free measures
(\ref{free}). The mixture is given in terms of the distribution of pinned
sites. For the $\delta$-pinning case, this is particularly simple.
$\mu_{\Lambda}^{\varepsilon}$ generates a law on subsets $A\subset\Lambda$,
the set of sites where the random field is $0$ inside $\Lambda.$ Conditioned
on this set, the field is then just the free field (\ref{free}) on
$A^{c}\overset{\mathrm{def}}{=}\Lambda\backslash A$ with $0$-boundary
conditions on $\left(  \mathbb{Z}^{d}\backslash\Lambda\right)  \cup A.$ It is
therefore crucial to have information on the distribution of pinned sites,
which we denote by $\nu_{\Lambda}^{\varepsilon}$ (see the precise definition
in (\ref{eq_exp_sw})). The main result on this problem is a domination
property of this distribution by Bernoulli measures from above and below. The
difficulty in dimension two (in contrast with the situation in higher
dimensions) is that, strictly speaking, there is no sharp domination, i.e. with
the same $\varepsilon$ behaviour from above and below,
but, surprisingly, correlations can be estimated as if there were such a
domination. This is the content of Theorem \ref{thm_obstacles} which is proved
for general convex interactions.

The main results on the depinning properties (Theorem \ref{thm_critical} and
Theorem \ref{thm_critical2}) are however proved for the Gaussian $\delta
$-pinning case only. The restriction to the Gaussian case is mainly due to the
fact that we need precise information on the behavior of various objects
appearing in the random walk representations (\ref{RW}), like estimates of Green
functions and ranges of the random walk. One might hope that with the help of
the Helffer-Sj\"{o}strand representation (see \cite{DGI}) which gives a
representation similar to (\ref{RW}) also for the case of convex, even
interactions, this could be extended. However, this random walk representation
is a much more complicated object and the precise information we need is not
available in this case, yet.

The restriction to the $\delta$-pinning case, which is made here mainly for
technical convenience, is more innocuous and could probably be much relaxed by
replacing the simple expansion of products by the more sophisticated
Brydges-Fr\"{o}hlich-Spencer random walk representation, see \cite{BFS} (not
to be mixed up with the Helffer-Sj\"{o}strand representation).

\medskip
The critical behavior of the 2-point correlation function has also been
obtained in a mean-field regime, mentioned at the beginning of the
introduction, in \cite{DMR}, see also \cite{Lemberger}. We briefly describe
the setting and the result in order to show the difference with the regime
studied here. The measure considered in~\cite{DMR} is
\[
\mu_{\Lambda}^{U}(\mathrm{d}\phi)\overset{\mathrm{def}}{=}\frac{1}{Z_{\Lambda
}^{U}}\;\prod_{\left\langle xy\right\rangle }e^{-\tfrac{1}{2}(\phi_{x}
-\phi_{y})^{2}}\,\prod_{x\in\Lambda}e^{-U(\phi_{x})}\;\prod_{x\in\Lambda
}\mathrm{d}\phi_{x}\prod_{y\not \in\Lambda}\delta_{0}(\mathrm{d}\phi_{y})\,,
\]
where $\left\langle xy\right\rangle $ denotes nearest-neighbor sites, and
\[
U(x)=-c(e^{-\frac{x^{2}}{2q^{2}}}-1)\,.
\]
Then, provided\footnote{It is emphasized in~\cite{DMR} that this condition is
actually too strong and that the result should be true under the weaker
condition that $K\log(1+c^{-1})<q$, which characterizes the mean-field
regime.} $K\log(1+c^{-1})<\sqrt{q}$ for some sufficiently large constant $K$
and $0<c\leq1$, it is proved that
\[
\mu^{U}(\phi_{x}\,\phi_{y})\leq K\log(q/\sqrt{c})e^{-D\,\frac{\sqrt{c}}
{q}\,\left|  x-y\right|  }\,,
\]
with the constant $D\rightarrow1$ if $c$ is fixed and $q\rightarrow\infty$.
The heuristics behind this result is rather clear. Under the above assumption,
the quadratic approximation $U(x)=\frac{c}{2q^{2}}\,x^{2}$ holds over a huge
range of values of $x$. Over this range of values the measure
$\mu^{U}$ behaves like a massive Gaussian model with mass $m=\sqrt{c}/q$, and
therefore, provided the interface stays mostly there, the
exponential decay should be given by this mass. The main part of the proof
in~\cite{DMR,Lemberger} was then to show perturbatively that indeed the
interface remains essentially all the time in this range.

The $\delta$-pinning corresponds to an opposite regime, where instead of
having a very wide and shallow potential well, one has a very narrow and deep
one. It is far less clear \emph{a priori} what the behavior of the correlation
lengths should be in this case, since the latter cannot be read from the
self-potential.

\medskip
The paper is organized as follows: In the next section we state precisely the
results. In Section \ref{SectDomination}, we prove the main domination
results. Section \ref{SectVariance} proves the results on the variance, and
\ref{SectMass} for the covariance. In Appendix~\ref{app_mass}, we prove the
existence of the mass in the Gaussian case. We will also need precise results
about standard random walks, and the number of points visited by random walks.
Some of these properties are standard, but others are more delicate. We
collect what we need in Appendix~\ref{app_RW} and Appendix \ref{RangeRW}. To
complete the picture, we shortly sketch the one-dimensional situation in
Appendix \ref{app_1d} but only in the nearest neighbor case $p(\pm1)=1/2$
which is easily reduced to standard renewal theory.

\textbf{Acknowledgments:} We thank Fran\c{c}ois Dunlop for several very
interesting discussions on the physical aspects of these questions, as well as
for encouraging us to look at the problem investigated in the present paper.
We also thank Pietro Caputo for interesting discussions. Y.V. gratefully
acknowledges the warm hospitality of the Institute for Mathematics of
Z\"{u}rich University where part of this work was done. Y.V. is supported by a
Swiss National Science Foundation Grant \#8220-056599. E.B. is supported by
SNSF Grant \#20-55648.98.

\section{Results \label{SectResults}}

The basic assumption is that the (symmetric) transition kernel $\left(
p(x)\right)  _{x\in\mathbb{Z}^{d}}$ is irreducible and satisfies
(\ref{CondH1}). Only in Theorem \ref{thm_critical2} we need a stronger
assumption. We write $X_{0},X_{1},X_{2},\ldots$ for a random walk with these
transition probabilities, and $\mathrm{P}_{x}$ for the corresponding law for a
walk starting in $x.$ With $X_{[0,n]}$ we denote the set of points visited by
the walk up to time $n,$ and by $\left|  X_{[0,n]}\right|  $ the number of
points visited. If $p(0)=0,$ then remark that the interface model is not
changed if we replace $p$ by its half, putting $p(0)=1/2,$ and doubling
$\beta.$ We can therefore as well assume that $p$ is aperiodic, and especially
that for any $x\in\mathbb{Z}^{2},$ $p_{n}(x)>0$ for large enough $n,$ where
$p_{n}$ is the $n$-fold convolution of $p.$

We denote by $C$ or $C^{\prime},C^{\prime\prime}$ generic constants, not
necessarily the same at different occurrences, which may depend on $p$ and the
dimension $d,$ but on nothing else, unless explicitly stated.

Our first result complements estimates obtained in \cite{BB, DV} where it was
shown for $d=2$ that provided $V^{\prime\prime}\geq\underline{c}>0$ and
$p(\,\cdot\,)$ satisfies (\ref{CondH1}), there exists a constant $C>0$
(depending on $p$ only) such that, for small enough $\underline{\mathfrak{e}%
}\overset{\mathrm{def}}{=}2a\sqrt{\beta\underline{c}}(e^{b}-1)>0$,
\[
\sup_{\Lambda}\mu_{\Lambda}^{a,b}(\phi_{0}^{2})\leq C\,\left(  \frac
{1}{\underline{c}\beta}\left|  \log\underline{\mathfrak{e}}\right|
+a^{2}\right)  \,.
\]
We are going to show that this upper bound indeed corresponds to the correct behavior:

\begin{theorem}
\label{thm_lowerbound}Assume $d=2$ and let $V$ be an even $C^{2}$ function
with $0\leq V^{\prime\prime}(x)\leq\overline{c}$ for all $x.$ Then in the
square well pinning case, there exists a constant $C>0$ (depending only on
$p)$ such that for any
\[
\overline{\mathfrak{e}}=\overline{\mathfrak{e}}\left(  a,b,\beta\right)
\overset{\mathrm{def}}{=}2a\sqrt{\beta\overline{c}}\left(  \mathrm{e}%
^{b}-1\right)
\]
small enough and provided $2a\sqrt{\beta\overline{c}}\leq|\log\overline
{\mathfrak{e}}|^{1/2}$,
\[
\liminf_{\Lambda\nearrow\mathbb{Z}^{2}}\mu_{\Lambda}^{a,b}\left(  \phi_{0}%
^{2}\right)  \geq\frac{C}{\overline{c}\beta}\left|  \log\overline{\mathfrak
{e}}\right|  .
\]
This remains true for $\delta$-pinning with $\overline{\mathfrak{e}%
}=\varepsilon\sqrt{\beta\overline{c}}$ (and $a=0$).
\end{theorem}

Our next two results are for the Gaussian (i.e. $V(x)=x^{2}/2)$ case and
$\delta$-pinning. For this case there is a simple proof of the existence of a
thermodynamic limit.

\begin{proposition}
\label{prop_thermolimit}The thermodynamic limit
\[
\mu^{\ast,\varepsilon}\overset{\mathrm{def}}{=}\lim_{\Lambda\nearrow
\mathbb{Z}^{2}}\mu_{\Lambda}^{\ast,\varepsilon}%
\]
exists in all dimensions and is translation invariant. The limit is defined in
terms of limits of integrals over bounded local functions.
\end{proposition}

\begin{proof}
This is an immediate consequence of the corresponding property for the law of
the pinned sites, given in Lemma \ref{LePinnedsites} below.
\end{proof}

Let $\mathcal{Q}$ be the covariance matrix of $p$: $\mathcal{Q}(i,j)\overset
{\mathrm{def}}{=}\sum_{x\in\mathbb{Z}^{d}}x_{i}x_{j}\,p(x).$ Our main result
on the behavior of the variance in the Gaussian case is the following

\begin{theorem}
\label{thm_critical} Assume $d=2.$ There exists $\varepsilon_{0}>0$ and $C>0 $
such that for all $\varepsilon$ and $\beta$ satisfying $0<\varepsilon
\sqrt{\beta}<\varepsilon_{0}$
\[
\left|  \mu^{\ast,\varepsilon}\left(  \phi_{0}^{2}\right)  -\frac{\left|
\log\left(  \sqrt{\beta}\varepsilon\right)  \right|  }{2\pi\beta\sqrt
{\det\mathcal{Q}}}\right|  \leq C\log\left|  \log\left(  \sqrt{\beta
}\varepsilon\right)  \right|
\]
\end{theorem}

The second quantity we are interested in is the decay-rate of the covariance
(i.e. the mass). This is of interest also in the higher-dimensional case. It
is defined, for $x$ on the unit sphere $x\in\mathbb{S}^{d-1}\overset
{\mathrm{def}}{=}\left\{  x\in\mathbb{R}^{d}:\left\|  x\right\|
_{2}=1\right\}  $ as the limit
\begin{equation}
m_{\varepsilon}(x)=-\lim_{k\rightarrow\infty}\frac{1}{k}\,\log\mu
^{\varepsilon}(\phi_{0}\,\phi_{\left[  kx\right]  }),\label{eq_mass}%
\end{equation}
where $\left[  kx\right]  $ is the integer part of $kx,$ componentwise. The
existence of this limit, in the Gaussian case, is proved in
Appendix~\ref{app_mass}. The following theorem shows that in the Gaussian case
$m_{\varepsilon}\sim\varepsilon^{1/2+o(1)}$ as $\varepsilon$ goes to zero,
provided the coupling $p(\,\cdot\,)$ has an exponential moment.

\begin{theorem}
\label{thm_critical2}Consider the case of $\delta$-pinning and Gaussian
interaction, and assume that there exists $a>0$ such that
\begin{equation}
\sum_{x\in\mathbb{Z}^{d}}p(x)\,\mathrm{e}^{a\left|  x\right|  }<\infty
\,.\label{CondH2}%
\end{equation}

a) Assume $d=2.$ Then there exist $\varepsilon_{0}>0$ and constants
$0<C_{1}\leq C_{2}<\infty$ (depending only on $p$) such that
\[
C_{1}\frac{\left(  \sqrt{\beta}\varepsilon\right)  ^{1/2}}{\left|  \log
(\sqrt{\beta}\varepsilon)\right|  ^{3/4}}\leq m_{\varepsilon}(x)\leq C_{2}%
\frac{\left(  \sqrt{\beta}\varepsilon\right)  ^{1/2}}{\left|  \log(\sqrt{\beta
}\varepsilon)\right|  ^{3/4}}%
\]
for all $0<\varepsilon\sqrt{\beta}<\varepsilon_{0}$ and for any $x\in
\mathbb{S}^{1}.$

b) Assume $d\geq3.$ Then there exist $\varepsilon_{0}>0$ and constants
$0<C_1<C_2<\infty$ (depending only on $p$ and $d$) such that
\[
C_1(\sqrt{\beta}\varepsilon)^{1/2}\leq m_{\varepsilon}(x)\leq
C_2(\sqrt{\beta}\varepsilon)^{1/2}%
\]
for all $0<\varepsilon\sqrt{\beta}<\beta_{0}$ and for any $x\in\mathbb{S}%
^{d-1}$.
\end{theorem}

\begin{remark}
1. The theorem gives much more than just the correct power law decay
$\varepsilon^{1/2}$ of the mass, since it shows that there is no logarithmic
corrections when $d\geq3$, while it provides the correct power for the
logarithmic correction when $d=2$. The most precise results one might expect to
hold in the latter case would be
\[
m_{\varepsilon}(x)=\frac{\left(  \sqrt{\beta}\varepsilon\right)  ^{1/2}%
}{\left|  \log\sqrt{\beta}\varepsilon\right|  ^{3/4}}\varphi\left(  x\right)
\left(  1+o(1)\right)  ,
\]
where $\varphi$ is a positive function on $\mathbb{S}^{1}$ which is bounded
and bounded away from $0$. Our techniques, however, do not give so precise an information.

2. The assumption on the existence of an exponential moment is essentially
optimal. Otherwise, there is no positive mass. Indeed, it is easy to show that
the decay of the covariance cannot be faster than that of $p(\,\cdot\,)$: In
the random-walk representation of $\mu^{\varepsilon}(\phi_{0}\,\phi_{x}) $,
see \eqref{eq_pincov}, we get a lower bound by letting the random-walk jump
directly from $0$ to $x$. Probably, this ``one-jump'' contribution gives the
leading order of the decay of correlations correctly, but we don't have a proof.
\end{remark}

\begin{remark}
\label{RemRescale}The temperature parameter enters only in a trivial way. If
we replace the field $\left(  \phi_{x}\right)  $ by $\left(  \sqrt{\beta}%
\phi_{x}\right)  ,$ and $\varepsilon$ by $\sqrt{\beta}\varepsilon$ we have
transformed the model to temperature parameter $\beta=1.$ In the proofs, we
will therefore always assume $\beta=1.$
\end{remark}

As remarked in the introduction, the mechanism at play is that the potential
will randomly pin some sites at height $0$ or close to $0$. The main point
therefore is to find the properties of the distribution of these pinned sites.
Precise information about this distribution is used in essential ways in the
proofs of the previous theorems. Since these results are also interesting
\textit{per se}, and yield a better understanding of the reason behind the
behavior described above, we discuss them in some details, and prove more than
is needed for the proofs of Theorems~\ref{thm_critical}
and~\ref{thm_critical2}. In particular, we do not restrict to the Gaussian case.

Let us start by defining precisely what we mean by the set of pinned sites,
and its distribution. The starting point is the following expansion: For any
bounded measurable function $f$,
\begin{align}
\mu_{\Lambda}^{a,b}(f)  &  =\frac{1}{Z_{\Lambda}^{a,b}}\int f(\phi
)\mathrm{e}^{-\tfrac{\beta}{2}\sum_{x,y}p(x-y)V(\phi_{x}-\phi_{y})}\prod
_{x\in\Lambda}\left\{  (\mathrm{e}^{b}-1)\mathrm{I}\left(  \left|  \phi
_{x}\right|  \leq a\right)  +1\right\}  \prod_{x\in\Lambda}\!\!\mathrm{d}%
\phi_{x}\prod_{x\not \in\Lambda}\!\delta_{0}(\mathrm{d}\phi_{x})\nonumber\\
&  =\sum_{A\subset\Lambda}(\mathrm{e}^{b}-1)^{\left|  A\right|  }%
\frac{Z_{\Lambda}^{a}(A)}{Z_{\Lambda}^{a,b}}\;\mu_{\Lambda}\left(  f\left|
\,\left|  \phi_{x}\right|  \leq a,\;\forall x\in A\right.  \right)
\label{eq_exp_sw}\\
&  =\sum_{A\subset\Lambda}\nu_{\Lambda}^{a,b}(A)\;\mu_{\Lambda}\left(
f\left|  \,\left|  \phi_{x}\right|  \leq a,\;\forall x\in A\right.  \right)
,\nonumber
\end{align}
where
\[
\nu_{\Lambda}^{a,b}(A)\overset{\mathrm{def}}{=}(\mathrm{e}^{b}-1)^{\left|
A\right|  }\frac{Z_{\Lambda}^{a}(A)}{Z_{\Lambda}^{a,b}},
\]%
\[
Z_{\Lambda}^{a}(A)\overset{\mathrm{def}}{=}Z_{\Lambda}\mu_{\Lambda}\left(
\left|  \phi_{x}\right|  \leq a,\,\forall x\in A\right)  .
\]
Therefore the effect of the potential can be seen as pinning, i.e.
constraining to the interval $[-a,a]$ a random set of points, the
\textit{pinned sites}. The distribution of the latter is given by the
probability measure $\nu_{\Lambda}^{a,b}$. We'll denote by $\mathcal{A}$ the
corresponding random variable, taking values in the subsets of $\Lambda$. A
completely similar representation is obtained in the case of $\delta$-pinning
by just expanding the term $\prod_{x\in\Lambda}(\mathrm{d}\phi_{x}%
+\varepsilon\delta_{0}(\mathrm{d}\phi_{x}))$. The result reads
\begin{equation}
\mu_{\Lambda}^{\varepsilon}(f)=\sum_{A\subset\Lambda}\nu_{\Lambda
}^{\varepsilon}(A)\;\mu_{A^{c}}(f)\;\,,\label{eq_exp_pinning}%
\end{equation}
where $A^{c}\overset{\mathrm{def}}{=}\Lambda\setminus A$ and $\nu_{\Lambda
}^{\varepsilon}(A)\overset{\mathrm{def}}{=}\varepsilon^{\left|  A\right|
}\frac{Z_{A^{c}}}{Z_{\Lambda}^{\varepsilon}}$.

The following lemma gives some basic properties of the distribution of pinned sites.

\begin{lemma}
\label{LePinnedsites} Suppose that Griffiths' inequalities (in the sense of~
\cite{Gi}) hold for the measure $\mu_{\Lambda}$. Then

1. $\nu_{\Lambda}^{a,b}$ and $\nu_{\Lambda}^{\varepsilon}$ satisfy the lattice
condition, i.e.
\begin{equation}
\nu_{\Lambda}^{a,b}\left(  A\cup B\right)  \nu_{\Lambda}^{a,b}\left(  A\cap
B\right)  \geq\nu_{\Lambda}^{a,b}\left(  A\right)  \nu_{\Lambda}^{a,b}\left(
B\right)  ,\label{FKG}%
\end{equation}
for $A,B\subset\Lambda,$ and similarly for $\nu_{\Lambda}^{\varepsilon}$. In
particular, these two measures are strong FKG, see~\cite{FKG}.

2. $\nu^{\varepsilon}\overset{\mathrm{def}}{=}\lim_{\Lambda\nearrow
\mathbb{Z}^{d}}\nu_{\Lambda}^{\varepsilon}$ exists and is translation invariant.
\end{lemma}

\begin{proof}
Part 1. is very simple: In the square-well case (\ref{FKG}) is equivalent to
\begin{align*}
&  \mu_{\Lambda}\left(  \left|  \phi_{x}\right|  \leq a,\,\forall x\in
B\setminus A\,\bigm\vert\,\left|  \phi_{y}\right|  \leq a,\,\forall y\in
A\right) \\
&  \geq\mu_{\Lambda}\left(  \left|  \phi_{x}\right|  \leq a,\,\forall x\in
B\setminus A\,\bigm\vert\,\left|  \phi_{y}\right|  \leq a,\,\forall y\in A\cap
B\right)
\end{align*}
which follows from Griffiths' inequality. The $\delta$-pinning case is similar.

Part 2 is easy, too: For any local increasing function $f$ (of the random set
$\mathcal{A}$) with support inside $\Lambda\subset\Lambda^{\prime}%
\Subset\mathbb{Z}^{d}$, one has
\begin{equation}
\nu_{\Lambda}^{\mathcal{\varepsilon}}(f)=\nu_{\Lambda^{\prime}}%
^{\mathcal{\varepsilon}}(f\,|\,\Lambda^{\prime}\setminus\Lambda\subset
\mathcal{A})\geq\nu_{\Lambda^{\prime}}^{\mathcal{\varepsilon}}%
(f)\,.\label{eq_monotone}%
\end{equation}
Translation invariance is a simple consequence of this. Indeed, let
$x\in\mathbb{Z}^{d}$ and $\mathcal{T}_{x}f=f(\,\cdot\,-x)$. Denoting by
$\Lambda_{+}$ (respectively $\Lambda_{-}$) the biggest (respectively smallest)
square box centered at $x$ contained in (respectively containing) $\Lambda$,
we have
\[
\nu_{\Lambda_{+}}^{\mathcal{\varepsilon}}(\mathcal{T}_{x}f)\geq\nu_{\Lambda
}^{\mathcal{\varepsilon}}(\mathcal{T}_{x}f)\geq\nu_{\Lambda_{-}}%
^{\mathcal{\varepsilon}}(\mathcal{T}_{x}f)\,,
\]
provided $\Lambda$ is big enough. Taking the limit $\Lambda\nearrow
\mathbb{Z}^{d}$ and using the fact that $\nu_{\Lambda_{-}}%
^{\mathcal{\varepsilon}}(\mathcal{T}_{x}f)=\nu_{\mathcal{T}_{-x}\Lambda_{-}%
}^{\mathcal{\varepsilon}}(f)$, and the corresponding statement for
$\Lambda_{+}$, we get $\nu^{\mathcal{\varepsilon}}(f)=\nu
^{\mathcal{\varepsilon}}(\mathcal{T}_{x}f)$ which implies the desired result.
\end{proof}

\begin{remark}
1. Griffiths' inequalities are known to hold in the Gaussian case,
see~\cite{Gi}.

2. Part 1. of the lemma is of course not specific to the cubic
lattice. Griffiths' inequality for $\mu_{\Lambda}$ implies the strong FKG
property for the distribution of pinned sites on an arbitrary lattice.
\end{remark}

The following Theorem \ref{thm_obstacles} is the key step for our analysis of
the random fields. It states domination properties of the field of pinned
sites by Bernoulli measures and is a substantial improvement on the results
already present in~ \cite{DV, IV}. Although the main emphasis in this paper is
on the case of the (difficult) two-dimensional lattice, we include also the
higher-dimensional case.

Let us first introduce some standard notions. If $\nu_{1}$ and $\nu_{2}$ are
two probability measures on the set of subsets $\left\{  0,1\right\}
^{\Lambda}$ of a finite set $\Lambda,$ we say that $\nu_{1}$
\textit{dominates} $\nu_{2},$ if for any increasing function $f:\mathcal{P}%
\left(  \Lambda\right)  \rightarrow\mathbb{R},$ we have
\begin{equation}
\nu_{1}\left(  f\right)  \geq\nu_{2}\left(  f\right)  .\label{Domination}%
\end{equation}
We say that $\nu_{1}$ \textit{strongly dominates} $\nu_{2},$ if for any
$x\in\Lambda$ and any subset $C\subset\Lambda\backslash\left\{  x\right\}  $%
\begin{equation}
\nu_{1}\left(  x\in\mathcal{A}\,|\,\mathcal{A\cap}\left(  \Lambda\setminus
\left\{  x\right\}  \right)  =C\right)  \geq\nu_{2}\left(  x\in\mathcal{A}%
\,|\,\mathcal{A\cap}\left(  \Lambda\backslash\left\{  x\right\}  \right)
=C\right)  .\label{StrongDomination}%
\end{equation}
It is evident that strong domination implies domination, and the latter
implies that for any subset $B\subset\Lambda,$ one has
\[
\nu_{1}\left(  \mathcal{A}\cap B=\emptyset\right)  \leq\nu_{2}\left(
\mathcal{A}\cap B=\emptyset\right)  .
\]

We formulate the next theorem for the square-well case only. We set
\begin{equation}
\varepsilon=\varepsilon_{a,b}\overset{\mathrm{def}}{=}2a\left(  \mathrm{e}%
^{b}-1\right)  .\label{DefEps}%
\end{equation}
The $\delta$-pinning case follows either in an identical way, or by taking the
limit as $a\rightarrow0,$ keeping $\varepsilon$ fixed.

\begin{theorem}
\label{thm_obstacles}Let $V$ be an even $C^{2}$ function.

1. Assume $d\geq2$ and suppose $0\leq V^{\prime\prime}(x)\leq\overline
{c},\,\forall x$. Then there exists $C<\infty$, depending only on $p$ and $d,$
such that for any $\Lambda\Subset\mathbb{Z}^{d}$, the distribution
$\nu_{\Lambda}^{a,b}$ of pinned sites is strongly dominated by the Bernoulli
measure on $\{0,1\}^{\Lambda}$ with density $p_{-}^{\prime}\overset
{\mathrm{def}}{=}C\;(1\wedge a^{-1})\sqrt{\beta\overline{c}}\,\varepsilon$
($\varepsilon$ given by (\ref{DefEps})). In particular, for any $B\subset
\Lambda\Subset\mathbb{Z}^{d}$,
\begin{equation}
\nu_{\Lambda}^{a,b}(\mathcal{A}\cap B=\emptyset)\geq(1-p_{-}^{\prime
})^{\left|  B\right|  }.\label{Bern1}%
\end{equation}

2. Assume $d=2$ and suppose that $V(x)=\tfrac{1}{2}x^{2}$. For any $\alpha>0
$, there exist $\varepsilon_{0}>0$ and $C\left(  \alpha\right)  <\infty$, such
that, for $\sqrt{\beta}\varepsilon\leq\varepsilon_{0},$ any $\Lambda
\Subset\mathbb{Z}^{2}$ and any $B\subset\Lambda\Subset\mathbb{Z}^{2}$ with
$d(B,\Lambda^{\mathrm{c}})>\varepsilon^{-\alpha}$
\begin{equation}
\nu_{\Lambda}^{a,b}(\mathcal{A}\cap B=\emptyset)\geq(1-p_{-})^{\left|
B\right|  },\label{Bern1b}%
\end{equation}
with
\begin{equation}
p_{-}=p_{-}\left(  \alpha,\varepsilon\right)  \overset{\mathrm{def}}%
{=}C\left(  \alpha\right)  \,|\log\sqrt{\beta}\varepsilon|^{-1/2}\sqrt{\beta
}\varepsilon.\label{PeMinus}%
\end{equation}

3. Assume $d=2$ and suppose $V^{\prime\prime}(x)\geq\underline{c}>0,\,\forall
x.$ There exist $\varepsilon_{0}>0$ and $C>0$ such that, for all $a,b>0$ with
$\sqrt{\beta\underline{c}}\,\varepsilon\leq\varepsilon_{0},\,2a\sqrt
{\beta\underline{c}}\leq|\log\sqrt{\beta\underline{c}}\,\varepsilon|^{1/2}$, and
for any set $B\subset\Lambda\Subset\mathbb{Z}^{2}$,
\begin{equation}
\nu_{\Lambda}^{a,b}(\mathcal{A}\cap B=\emptyset)\leq(1-p_{+})^{\left|
B\right|  }\,,\label{Bern2}%
\end{equation}
with
\begin{equation}
p_{+}\overset{\mathrm{def}}{=}C\left|  \log\sqrt{\beta\underline{c}%
}\,\varepsilon\right|  ^{-1/2}\sqrt{\beta\underline{c}}\,\varepsilon
.\label{PePlus}%
\end{equation}

4. For $d\geq3$ and $V^{\prime\prime}(x)\geq\underline{c}>0,\,$there exists
$C>0,$ depending only on $p$ and $d$ such that $\nu_{\Lambda}^{a,b}$ strongly
dominates a Bernoulli measure with
\begin{equation}
p_{+}\overset{\mathrm{def}}{=}C\left(  1\wedge a^{-1}\right)  \sqrt
{\beta\underline{c}}\,\varepsilon\label{PePlus3}%
\end{equation}

All the statements remain true in the case of $\delta$-pinning.
\end{theorem}

\begin{remark}
1. Part 3 of the theorem is stated for small enough $\varepsilon$ and $a$ only.
An essentially identical proof yields exponential decay of $\nu_{\Lambda
}^{a,b}(\mathcal{A}\cap B=\emptyset)$ for any $a,b>0$. The precise
$\varepsilon$ dependence given in the theorem, however, is only valid for
small values of $\varepsilon$.\\
2. We expect that part 2 could be generalized to more general convex
interactions $V$, but a proof eludes us.
\end{remark}

The fact that for $d\geq3,$ $\nu_{\Lambda}^{a,b}$ can be strongly dominated
from above \textit{and below} by a Bernoulli measure has been observed by Dima
Ioffe (oral communication). That this is not true for $d=2$ can be seen as
follows: It is easy to check that
\[
\nu_{\Lambda}^{\varepsilon}\left(  \mathcal{A}\ni x\,\left|  y\not \in
\mathcal{A},\,\forall y\neq x\text{ s.t. }\left|  x-y\right|  <T\right.
\right)
\]
is decreasing to zero as $T\rightarrow\infty,$ $\Lambda\nearrow\mathbb{Z}^{2}%
$, since under this conditioning typical values of the field at the sites
neighboring $x$ will be (at least) of order $\sqrt{\log T}$. This excludes the
possibility of \emph{any} strong domination of a Bernoulli measure, uniformly
in $\Lambda$. This leaves open the possibility of a domination in the sense of
(\ref{Domination}), which might be true; note however that the density of the
corresponding Bernoulli measure cannot be larger than
$\varepsilon|\log\varepsilon|^{-1/2}$.

The difference between $p_{-}^{\prime}$ and $p_{-}$ as function of
$\varepsilon$ for $d=2$ in part 1 and 2 is necessary. Indeed, there is no
strong domination by a Bernoulli process with density $o(\varepsilon)$, as the
following argument shows: Consider $\delta$-pinning; it is enough to show
that
\[
\nu_{\Lambda}^{a,b}(\mathcal{A}\ni0\,|\,\Lambda\setminus\{0\}\subset
\mathcal{A})\geq C\varepsilon\,,
\]
but this is obvious since
\[
\nu_{\Lambda}(\mathcal{A}\ni0\,|\,\Lambda\setminus\{0\}\subset\mathcal{A}%
)=\left(  1+\varepsilon^{-1}\mathrm{Z}_{\{0\}}\right)  ^{-1}\,.
\]
In fact, even more is true: There is no domination, even in the sense of
(\ref{Domination}), by a Bernoulli measure of density $o(\varepsilon)$. Indeed,
it is not difficult to show that the indicator function of the increasing event
$\{\mathcal{A}\supset B\}$ is larger than
$$
(C\varepsilon)^{|B|}\,|\log\varepsilon|^{-1}\,,
$$
for any connected set $B\subset\Lambda$. This shows in particular that there
must be a gap between any upper and lower domination of
$\nu_\Lambda^\varepsilon$ in dimension $2$. In view of this, it is rather
remarkable that as long as we are only interested in covariances of the field,
such a domination holds, as a consequence of the estimates (\ref{Bern1b}) and
(\ref{Bern2}):

\begin{corollary}
\label{CorDom} Assume the Gaussian $\delta$-pinning case with $\beta=1$ (which
is no restriction, according to Remark \ref{RemRescale}). There exists
$\varepsilon_{0}>0$ such that for $0<\varepsilon\leq\varepsilon_{0}$ the
following is true. Let $\rho_{+}$ be the Bernoulli measure with density
(\ref{PePlus}) or (\ref{PePlus3}), and $\rho_{-}$ the Bernoulli measures with
density $p_{-}\left(  \varepsilon\right)  $ from (\ref{PeMinus}) in the case
$d=2$, and $p_{-}^{\prime}\left(  \varepsilon\right)  $ in
the case $d\geq3$ $,$ then for any $x,y\in\mathbb{Z}^{d}$,
\[
\mu^{\star,\varepsilon}(\phi_{x}\,\phi_{y})\geq\rho_{-}\left(
\mu_{\mathcal{A}^{c}}^{\star}(\phi_{x}\,\phi_{y})\right)  \,,
\]
and
\[
\mu^{\star,\varepsilon}(\phi_{x}\,\phi_{y})\leq\rho_{+}\left(
\mu_{\mathcal{A}^{c}}^{\star}(\phi_{x}\,\phi_{y})\right)  \,.
\]
\end{corollary}

\begin{proof}
We recall that the variance of the Gaussian field can be written
\begin{equation}
\mu_{\Lambda}^{\star}(\phi_{x}\,\phi_{y})=\beta^{-1}\,\sum_{n\geq0}%
\mathrm{P}_{x}[X_{n}=y,\tau_{\Lambda}>n]\,,\label{eq_gausscov}%
\end{equation}
where $P_{x}$ is the law of the random walk in $\mathbb{Z}^{d}$, with
transition probabilities $p(\,\cdot\,)$, starting at $x$. Inserting this
in~\eqref{eq_exp_pinning}, we get
\begin{equation}
\mu_{\Lambda}^{\star,\varepsilon}(\phi_{x}\,\phi_{y})=\beta^{-1}\,\sum
_{n\geq0}\nu_{\Lambda}^{\varepsilon}\otimes\mathrm{P}_{x}[X_{n}=y,\tau
_{\mathcal{A}^{c}}>n]\,.\label{eq_pincov}%
\end{equation}
(Remember that $D^{c}=\Lambda\setminus D$.) Taking the expectation w.r.t.
$\nu_{\Lambda}^{\varepsilon}$ inside, we get
\[
\mu_{\Lambda}^{\star,\varepsilon}(\phi_{x}\,\phi_{y})=\beta^{-1}\,\sum
_{n\geq0}\mathrm{E}_{x}\left[  \mathrm{I}\left(  X_{n}=y\right)
\,\mathrm{I}\left(  \tau_{\Lambda}>n\right)  \,\nu_{\Lambda}^{\varepsilon
}(\mathcal{A}\cap X_{[0,n]}=\emptyset)\right]  \,.
\]
The corollary then follows from an application of the estimates of
Theorem~\ref{thm_obstacles}.
\end{proof}

Notice that Corollary \ref{CorDom} can also be stated in the two following
ways:
\begin{equation}
\sum_{n\geq0}\rho_{-}\otimes\mathrm{P}_{x}[X_{n}=x,T_{\mathcal{A}}>n]\leq
\mu^{\star,\varepsilon}(\phi_{x}\,\phi_{y})\leq\sum_{n\geq0}\rho_{+}\otimes
\mathrm{P}_{x}[X_{n}=y,\,T_{\mathcal{A}}>n],\label{eq_bernoulliRE}%
\end{equation}
where $T_{B}\overset{\mathrm{def}}{=}\min\left\{  n\geq0:X_{n}\in B\right\}
$, and, setting $\widetilde{p}_-$ equal to $p_-$ when $d=2$ and $p'_-$ when
$d\geq3$,
\begin{equation}
\sum_{n\geq0}\mathrm{E}_{x}\left[  \mathrm{I}\left(  X_{n}=y\right)  \,\left(
1-\widetilde{p}_{-}\right)  ^{\left|  X_{\left[  0,n\right]  }\right|  }\right]  \leq
\mu^{\star,\varepsilon}(\phi_{x}\,\phi_{y})\leq\sum_{n\geq0}\mathrm{E}_{x}\left[
\mathrm{I}\left(  X_{n}=y\right)  \,(1-p_{+})^{\left|  X_{\left[  0,n\right]
}\right|  }\right]  \,.\label{eq_GFsausage}%
\end{equation}

The problem is therefore essentially reduced to the analysis of the
asymptotics of the Green function of the random walk
with transition probabilities $p(\,\cdot\,)$, in an annealed random
environment of killing obstacles distributed according to Bernoulli measures
in the limit of vanishing density. Equivalently, what we need is the asymptotics of
the Green function of the ``Wiener sausage'',
\[
\sum_{n\geq0}\mathrm{E}_{x}\left[  \mathrm{I}\left(  X_{n}=y\right)
\,\mathrm{e}^{-s\left|  X_{\left[  0,n\right]  }\right|  }\right]  \,,
\]
as $s\rightarrow0$.

Let us conclude by making some comments on open problems. First of all, one
might wonder how universal the asymptotic behavior we have found actually is.
It would be very interesting to extend the analysis to a more general class of
interactions $V$. As remarked in the introduction, for even, strictly convex,
$C^{2}$ interactions a representation of the covariance, similar
to~\eqref{eq_gausscov}, also exists~\cite{DGI}. It was used in particular to
establish exponential decay of covariances for this class of
interactions~\cite{IV}. It is however much more complicated than the standard
random walk: The jump-rates of the walk are random, both in space and time;
they are given by the state of an independent diffusion process on
$\mathbb{R}^{\mathbb{Z}^{d}}$ which depends on the distribution of pinned
sites. So, even though the distribution of pinned sites can be treated in
general (see Theorem~\ref{thm_obstacles}), precise asymptotics in this
situation are probably hard to obtain.

Finally, there is a natural extension of this problem, which is more closely
related to the issue of critical wetting discussed in the beginning of the
paper: what happens in the presence of a hard-wall condition? More precisely,
one considers the measure
\[
\mu_{\Lambda}^{a,b,+}\overset{\mathrm{def}}{=}\mu_{\Lambda}^{a,b}%
(\,\cdot\,|\,\phi_{x}\geq0,\,\forall x\in\mathbb{Z}^{d})\,,
\]
or the corresponding measure with $\delta$-pinning. In this case, attraction
of the pinning potential competes with entropic repulsion due to the
conditioning, which makes this a much more difficult problem. Up to now, the
only rigorous results (in dimension larger than $1$) concern the existence, or
not, of a strictly positive critical value $\varepsilon_{c}$ such that for
$\varepsilon>\varepsilon_{c}$ the interface is pinned, while it is repelled
for $0<\varepsilon<\varepsilon_{c}$. It was shown in~\cite{BDZ} that for
quadratic interactions and dimensions $3$ and higher, there is no such
$\varepsilon_{c}$: As in the pure pinning case, the interface is localized for
arbitrarily weak pinning strength. On the other hand, it was shown
in~\cite{CV} that in dimension $2$ there exists such an $\varepsilon_{c}$;
moreover it was shown in the latter paper that this is true in any dimension
if the interaction is Lipschitz. The results of these two papers provide only
information on the density of pinned sites, but give no local estimates. For
example, it is even an open problem whether in the localized regime the
variance of the spin at the origin is finite. To get much more, namely the
critical behavior of such a quantity, seems therefore to be quite a challenge.

\section{Geometry of the pinned sites: Proof of Theorem~\ref{thm_obstacles}
\label{SectDomination}}

Note that it is enough to consider the case $\beta=1$ and $V^{\prime\prime
}\leq1$, respectively $V^{\prime\prime}\geq1$, in point 1, respectively 3 and 4.
Indeed, say in point 1, we can define $\widetilde{V}(x)=\beta V(x\big/%
\sqrt{\beta\overline{c}})$, and then, by an obvious change of variables we see
that
\begin{equation}
\nu_{\Lambda,\beta,V}^{a,b}=\nu_{\Lambda,1,\widetilde{V}}^{a\sqrt
{\beta\overline{c}},b}\label{EqRescale}%
\end{equation}
and by construction $\widetilde{V}^{\prime\prime}\leq1$.

\subsection{Proof of point 1.}

By simple algebraic manipulations, one can write, for any $A\subset
\Lambda\setminus\{x\}$,
\begin{equation}
\nu_{\Lambda}^{a,b}(\mathcal{A}\not \ni x\,|\,\mathcal{A}=A\text{ off
}x)=\left\{  1+(\mathrm{e}^{b}-1)\frac{\mathrm{Z}_{\Lambda}(A\cup
\{x\})}{\mathrm{Z}_{\Lambda}(A)}\right\}  ^{-1}\,.\label{eq_rem}%
\end{equation}
We now need the following result, which we establish below,
\begin{equation}
\frac{\mathrm{Z}_{\Lambda}(A\cup\{x\})}{\mathrm{Z}_{\Lambda}(A)}\leq
2a\;\frac{\mathrm{Z}_{\Lambda\setminus\{x\}}(A)}{\mathrm{Z}_{\Lambda}(A)}%
\leq2a\;\sqrt{\frac{1}{2\pi}}\,.\label{eq_density}%
\end{equation}
Of course, we also have the trivial upper bound $\mathrm{Z}_{\Lambda}%
(A\cup\{x\})/\mathrm{Z}_{\Lambda}(A)\leq1$, since the ratio can be written as
a conditional probability. This and (\ref{eq_density}) readily imply the
claim, since
\[
\nu_{\Lambda}^{a,b}(\mathcal{A}\not \ni x\,|\,\mathcal{A}=A\text{ off }%
x)\geq\bigl(1+C(1\wedge a^{-1})\varepsilon\bigr)^{-1}\geq1-C(1\wedge
a^{-1})\varepsilon\,.
\]
Let us now prove (\ref{eq_density}). The first inequality follows from the
fact that the maximum of the density $F_{\Lambda,A}$ of $\phi_{x}$ under
$\mu_{\Lambda}(\,\cdot\,|\,\left|  \phi_{z}\right|  \leq a,\,\forall z\in A) $
is at $\phi_{x}=0$. Indeed, $F_{\Lambda,A}^{\prime}(t)$ is equal to
\[
C_{\Lambda}(A,t)\,\sum_{y\in\Lambda}\,p(y-x)\,\mu_{\Lambda}\bigl(V^{\prime
}(\phi_{y}-t)\,\bigm\vert\,\phi_{x}=t,\,\left|  \phi_{z}\right|  \leq
a,\,\forall z\in A\bigr)-\sum_{y\not \in\Lambda}\,p(y-x)\,V^{\prime}(t)\,,
\]
where $C_{\Lambda}(A,t)>0$. Now, $V^{\prime}(s)\geq0$ for all $s\geq0$, and,
for $t\geq0$,
\begin{align*}
\mu_{\Lambda}\bigl(V^{\prime}(\phi_{y}-t)\,\bigm\vert\,\phi_{x}=t,\,\left|
\phi_{z}\right|  \leq a,\,\forall z\in A\bigr) &  =\mu_{\Lambda}^{-t}\bigl
(V^{\prime}(\phi_{y})\,\bigm\vert\,\phi_{x}=t,\,\left|  \phi_{z-t}\right|
\leq a,\,\forall z\in A\bigr)\\
&  \leq\mu_{\Lambda}\bigl(V^{\prime}(\phi_{y})\,\bigm\vert\,\phi
_{x}=0,\,\left|  \phi_{z}\right|  \leq(a-t)\vee0,\,\forall z\in A\bigr)\\
&  =0\,,
\end{align*}
where $\mu_{\Lambda}^{-t}$ denotes the measure with boundary condition $-t$
outside $\Lambda.$ The inequality is a consequence of FKG property, and the
last equality follows from the fact that $V^{\prime}$ is odd. Since
$F_{\Lambda,A}$ is even, the claim is proven.

To prove the second inequality in~\eqref{eq_density}, we write
\begin{align*}
\frac{\mathrm{Z}_{\Lambda}(A)}{\mathrm{Z}_{\Lambda\setminus\{x\}}(A)}  &
=\mu_{\Lambda\setminus\{x\}}\bigl(\int_{-\infty}^{\infty}\exp\bigl[-\sum
_{y\in\mathbb{Z}^{2}}p(y-x)(V(\phi_{y}-t)-V(\phi_{y}))\bigr]\;\mathrm{d}
t\,\Big\vert\,\left|  \phi_{z}\right|  \leq a,\,\forall z\in A\bigr)\\
&  \geq\int_{-\infty}^{\infty}\exp\bigl[-\tfrac{1}{2}\sum_{y\in\mathbb{Z}^{2}
}p(y-x)\\
&  \hspace{2cm}\times\mu_{\Lambda\setminus\{x\}}\left(  V(\phi_{y}
-t)+V(\phi_{y}+t)-2V(\phi_{y})\Bigm\vert\left|  \phi_{z}\right|  \leq
a,\,\forall z\in A\right)  \bigr]\mathrm{d}t\\
&  \geq\int_{-\infty}^{\infty}\exp[-\tfrac{1}{2}\,t^{2}]\,\mathrm{d}t,
\end{align*}
where the first inequality is a consequence of Jensen's inequality and the
symmetry of the measure under $\phi\rightarrow-\phi$, and for the second
inequality we used the assumption $V^{\prime\prime}\leq1$.

\subsection{Proof of point 2.}

We assume $d=2$ in this subsection. Let's write $B=\{t_{1},\dots,t_{|B|}\}$,
and let $B_{0}\overset{\mathrm{def}}{=}\emptyset$, $B_{k}=\{t_{1},\dots
,t_{k}\}$. Let also $C_{k}=\{x\in\Lambda\,|\,|x-t_{k}|\leq\varepsilon
^{-(\alpha\wedge\tfrac{1}{3})}\}$. We write
\begin{align*}
\nu_{\Lambda}\left(  \mathcal{A}\cap B=\emptyset\right)   &  =\prod
_{k=1}^{|B|}\nu_{\Lambda}\left(  \mathcal{A}\cap B_{k}=\emptyset
\,|\,\mathcal{A}\cap B_{k-1}=\emptyset\right) \\
&  =\prod_{k=1}^{|B|}\nu_{\Lambda}\left(  \mathcal{A}\not \ni t_{k}%
\,|\,\mathcal{A}\cap B_{k-1}=\emptyset\right)  \,.
\end{align*}
Now,
\begin{align*}
\nu_{\Lambda}\left(  \mathcal{A}\not \ni t_{k}\,|\,\mathcal{A}\cap
B_{k-1}=\emptyset\right)   &  =\Biggl\{1+(\mathrm{e}^{b}-1)\frac
{\sum_{A\not \ni t_{k},\,A\cap B_{k-1}=\emptyset}(\mathrm{e}^{b}%
-1)^{|A|}Z_{\Lambda}(A\cup\{t_{k}\})}{\sum_{A\not \ni t_{k},\,A\cap
B_{k-1}=\emptyset}(\mathrm{e}^{b}-1)^{|A|}Z_{\Lambda}(A)}\Biggr\}^{-1}\\
&  =\Biggl\{1+(\mathrm{e}^{b}-1)\frac{\nu_{\Lambda}\left(  \mathrm{I}%
(\mathcal{A}\not \ni t_{k})\,\frac{Z_{\Lambda}(\mathcal{A}\cup\{t_{k}%
\})}{Z_{\Lambda}(\mathcal{A})}\,\bigm|\,\mathcal{A}\cap B_{k-1}=\emptyset
\right)  }{\nu_{\Lambda}\left(  \mathcal{A}\not \ni t_{k}\,\bigm|%
\,\mathcal{A}\cap B_{k-1}=\emptyset\right)  }\Biggr\}^{-1}\,.
\end{align*}
Strong domination by Bernoulli measure from part 1 of the theorem shows that
$$
\nu_{\Lambda}\left(  \mathcal{A}\not \ni t_{k}\,|\,\mathcal{A}\cap
B_{k-1}=\emptyset\right)  \geq1/2\,,
$$
provided $\varepsilon$ is small enough. We
are left with the numerator. We decompose it as follows:
\begin{multline}
\nu_{\Lambda}\Bigl(  \mathrm{I}(\mathcal{A}\not \ni t_{k})\,\frac
{Z_{\Lambda}(\mathcal{A}\cup\{t_{k}\})}{Z_{\Lambda}(\mathcal{A})}\,\Bigm
|\,\mathcal{A}\cap B_{k-1}=\emptyset\Bigr)\\
=\nu_{\Lambda}\Bigl(  \mathrm{I}(\mathcal{A}\cap C_{k}=\emptyset
)\,\frac{Z_{\Lambda}(\mathcal{A}\cup\{t_{k}\})}{Z_{\Lambda}(\mathcal{A}%
)}\,\Bigm|\,\mathcal{A}\cap B_{k-1}=\emptyset\Bigr) \label{eq_decomp}\\
+\nu_{\Lambda}\Bigl(  \mathrm{I}(\mathcal{A}\cap C_{k}\neq\emptyset
,\,\mathcal{A}\not \ni t_{k})\,\frac{Z_{\Lambda}(\mathcal{A}\cup\{t_{k}%
\})}{Z_{\Lambda}(\mathcal{A})}\,\Bigm|\,\mathcal{A}\cap B_{k-1}=\emptyset
\Bigr)
\end{multline}
Let us first consider the second term. We already know, see~(3.3), that
$Z_{\Lambda}(A\cup\{t_{k}\})/Z_{\Lambda}(A)\leq2a/\sqrt{2\pi}$, for all
$A\not \ni t_{k}$. Therefore applying again the domination result from part 1,
this term is bounded from above by
\begin{align*}
\frac{2a}{\sqrt{2\pi}}\nu_{\Lambda}\left(  \mathcal{A}\cap C_{k}\neq
\emptyset\,\bigm|\,\mathcal{A}\cap B_{k-1}=\emptyset\right)   & =\frac
{2a}{\sqrt{2\pi}}\left(  1-\nu_{\Lambda}\left(  \mathcal{A}\cap C_{k}%
=\emptyset\,\bigm|\,\mathcal{A}\cap B_{k-1}=\emptyset\right)  \right) \\
& \leq\frac{2a}{\sqrt{2\pi}}\left(  1-(1-p_{-})^{|C_{k}|}\right)
=C\,2a\,\varepsilon^{1/3}.
\end{align*}
Let us now examine the first term in (\ref{eq_decomp}). We prove below that
\begin{align}
\frac{Z_{\Lambda}(A\cup\{t_{k}\})}{Z_{\Lambda}(A)}  & =\mu_{\Lambda}\left(
|\phi_{t_{k}}|\leq a\,|\,|\phi_{x}|\leq a,\,\forall x\in A\right)
\label{eq_byFKG}\\
& \leq\mu_{A^{\mathrm{c}}}\left(  |\phi_{t_{k}}|\leq2a\right)  \,.\nonumber
\end{align}
This then implies the following bound
\[
\frac{Z_{\Lambda}(A\cup\{t_{k}\})}{Z_{\Lambda}(A)}\leq C\,2a\,|\log
\varepsilon|^{-1/2}\,,
\]
since\footnote{This is the only place where we use the assumption that V is
quadratic. We don't know how to estimate the probability density in the
non-Gaussian case. Note that we only need to estimate it at zero, since the
maximum is there.}, under $\mu_{A^{\mathrm{c}}}$ with $A\cap C_{k}=\emptyset$,
$\phi_{t_{k}}$ is a Gaussian random variable with $0$ mean and variance
bounded from below by $C|\log\varepsilon|$. Putting all this together, we get
\[
\nu_{\Lambda}\left(  \mathcal{A}\not \ni t_{k}\,\bigm|\,\mathcal{A}\cap
B_{k-1}=\emptyset\right)  \geq\left\{  1+C\,|\log\varepsilon|^{-1/2}%
\,\varepsilon+C\,\varepsilon^{4/3}\right\}  ^{-1}\geq e^{-C\,|\log
\varepsilon|^{-1/2}\,\varepsilon}\,,
\]
and therefore
\[
\nu_{\Lambda}\left(  \mathcal{A}\cap B=\emptyset\right)  \geq e^{-C\,|\log
\varepsilon|^{-1/2}\,\varepsilon\,|B|}\geq\bigl(1-C\,|\log\varepsilon
|^{-1/2}\,\varepsilon\bigr)^{|B|}\,.
\]
It only remains to prove~\eqref{eq_byFKG}.
\[
\mu_{\Lambda}\left(  |\phi_{t_{k}}|\leq a\,\bigm|\,|\phi_{x}|\leq a,\,\forall
x\in A\right)  =1-2\mu_{\Lambda}\left(  \phi_{t_{k}}>a\,\bigm|\,|\phi_{x}|\leq
a,\,\forall x\in A\right)
\]
We use the FKG inequality, stating that the random field $\left(  \phi
_{x}\right)  _{x\in\Lambda\backslash A}$ with boundary conditions $\left\{
\phi_{x}=\eta_{x}:x\in A\cup\Lambda^{c}\right\}  ,$ depends monotonely on
$\left(  \eta_{x}\right)  _{x\in A\cup\Lambda^{c}}.$ Therefore, for $\left|
\eta_{x}\right|  \leq a,$ $x\in A$%
\begin{align*}
\mu_{\Lambda}\left(  \phi_{t_{k}}>a\,\bigm|\,\phi_{x}=\eta_{x},\,\forall x\in
A\right)   &  \geq\mu_{\infty}\left(  \phi_{t_{k}}>a\,\bigm|\,\phi
_{x}=-a,\,\forall x\in A\cup\Lambda^{c}\right) \\
&  =\mu_{A^{c}}\left(  \phi_{t_{k}}>2a\right)  =\frac{1-\mu_{A^{c}}\left(
\left|  \phi_{t_{k}}\right|  \leq2a\right)  }{2}.
\end{align*}
This proves (\ref{eq_byFKG}).

\subsection{Proof of point 3.}

\label{ssec_obstacles_2} We again have the assumption $d=2.$ The proof
proceeds in three steps. First, we prove a statement similar to that of
Theorem~\ref{thm_obstacles}, but valid only for sets $B$ sufficiently ``fat''.
In the second step, we use this result to show that with high probability
there is a high density of pinned sites at a large enough
($\mathcal{\varepsilon}$-dependent) scale. Then, in the last step, we use this
information to conclude the proof of Theorem~\ref{thm_obstacles}, part 3.

We need the following definition: Consider a partition of $\mathbb{Z}^{2}$
into cells by a grid of spacing $l$; the set of all cells entirely contained
in a subset (not necessarily finite) $\Lambda\subset\mathbb{Z}^{2}$ is denoted
by $\Lambda(l)$.

\subsubsection{Step 1: Probability of clean fat sets}

This step is a variant of the proofs given in~\cite{DV,IV}. Here, however, we
want to keep track of the $\mathcal{\varepsilon}$-dependence of the constants.
We remind the reader that we assume $\beta=1,$ \underline{$c$}$=1, $ and that
$\varepsilon=\varepsilon_{a,b}=2a\left(  \mathrm{e}^{b}-1\right)  .$

\begin{proposition}
\label{pro_fat} Let $\beta=1$ and let $V$ be an even, $C^{2}$, function with
$V^{\prime\prime}(x)\geq1$. There exist a constants $K>0$ and $\varepsilon
_{0}>0$ such that, for all $\varepsilon\leq\varepsilon_{0}$, and provided
$2a\leq|\log\varepsilon|$, the following holds. For any set $B\subset
\Lambda\Subset\mathbb{Z}^{2}$ composed of cells of $\Lambda(K${$\left|
\log\varepsilon\right|  $}${^{1/4}\,\varepsilon^{-1/2}})$
\[
\nu_{\Lambda}^{a,b}(\mathcal{A}\cap B=\emptyset)\leq\exp\left[  -C\left|
\log\varepsilon\right|  ^{-1/2}\,\varepsilon\left|  {B}\right|  \right]  \,.
\]
This statement remains true in the case of $\delta$-pinning.
\end{proposition}

\begin{proof}
We suppose first, for simplicity, that $B$ is connected. The changes for the
general case are the same as those described in~\cite{IV}, and we'll indicate
their effects on our bounds at the end of the proof.

Let $B^{0}\overset{\mathrm{def}}{=}B$, and define $B^{k+1}$ as the union of
$B^{k}$ and all its nearest neighboring cells in $\mathbb{Z}^{2}(K\left|
\log\varepsilon\right|  {^{1/4}\,\varepsilon^{-1/2}})$; let $\overline{k}$ be
the largest $k$ for which $B^{k}\subset\Lambda$. We then write
\[
\nu_{\Lambda}^{a,b}(\mathcal{A}\cap B=\emptyset)\leq\sum_{k=0}^{\overline{k}%
}\nu_{\Lambda}^{a,b}(\mathcal{A}\cap B^{k}=\emptyset\,|\,\mathcal{A}\cap
B^{k+1}\neq\emptyset)\,,
\]
and
\begin{align}
\nu_{\Lambda}^{a,b}(\mathcal{A}\cap B^{k}  & =\emptyset\,|\,\mathcal{A}\cap
B^{k+1}\neq\emptyset)\label{eq_probclean}\\
& \leq\left\{  \sum_{D\subset B^{k}}(\mathrm{e}^{b}-1)^{\left|  D\right|
}\inf\nolimits_{\substack{A\cap B^{k}=\emptyset\\ A\cap B^{k+1}\neq\emptyset
}}\mu_{\Lambda}\left(  \left|  \phi_{x}\right|  \leq a,\,\forall x\in
D\,\bigm\vert\,\left|  \phi_{x}\right|  \leq a,\,\forall x\in A\right)
\right\}  ^{-1}.\nonumber
\end{align}
It was proved in~\cite{DV}, see the proof of Proposition~4.1, that
\begin{align*}
\inf_{\substack{A\cap B^{k}=\emptyset\\ A\cap B^{k+1}\neq\emptyset}%
}\mu_{\Lambda}^{a,b}\left(  \left|  \phi_{x}\right|  \leq a,\,\forall x\in
D\,\bigm\vert\,\left|  \phi_{x}\right|  \leq a,\,\forall x\in A\right)   &
\geq\left(  C\left(  \frac{2a}{\sqrt{\left|  \log\varepsilon\right|  }}%
\wedge1\right)  \right)  ^{\left|  D\right|  }\\
&  =\left(  C\frac{2a}{\sqrt{\left|  \log\varepsilon\right|  }}\right)
^{\left|  D\right|  }\,,
\end{align*}
for the class of sets $D$ containing exactly one point in each cell of $B^{k}
$. Therefore, summing only over such $D$'s in (\ref{eq_probclean}) (notice
that there are $K^{2}${$\varepsilon$}$^{-1}\left|  \log\varepsilon\right|
^{1/2}$ choices for which site is occupied in a given cell), we get, choosing
$K^{2}=2/C,$ ($C$ from the formula above),
\[
\nu_{\Lambda}^{a,b}(\mathcal{A}\cap B^{k}=\emptyset\,|\,\mathcal{A}\cap
B^{k+1}\neq\emptyset)\leq\exp\left[  -C^{\prime}\,\left|  \log\varepsilon
\right|  ^{-1/2}\,\varepsilon\,\left|  B^{k}\right|  \right]  .
\]
From this we easily prove the claim for the one-component case, by summing
over $k$. Indeed, we can use the trivial estimate $\left|  B^{k}\right|
\geq\left|  B\right|  +k\varepsilon^{-1}\left|  \log\varepsilon\right|  ^{1/2}.$

To treat the case of multiple components, one proceeds as in the proof of
Theorem~2 in~\cite{IV}. The idea is to grow simultaneously all components in a
suitable way. This procedure only modifies the value of the constant in the
exponent, provided the components are all big enough. In our present
situation, this is enforced automatically as soon as {$\mathcal{\varepsilon}
$} is sufficiently small (the cells from which $B$ is built are growing when
{$\mathcal{\varepsilon}$} decreases).
\end{proof}

\subsubsection{Step 2: Density of pinned sites at large scales}

Our aim in this step is to show that any subset of $\Lambda$ has the property
that many of its points are close to pinned sites. To do this, we need two
partitions of $\mathbb{Z}^{2}$, first the one used in Step 1, $\mathbb{Z}%
^{2}(K\left|  \log\varepsilon\right|  {^{1/4}\,}${$\varepsilon$}${^{-1/2}})$,
and a second $\mathbb{Z}^{2}(\left|  \log\varepsilon\right|  \varepsilon
^{-1/2})$. The cells of the latter are called ``big'', and are supposed to be
built of cells from the finer partition (this might require some slight
modification of the size of the cells, but this is a trivial point). The
actual choice of the size of the big cells is actually not important. $\left|
\log\varepsilon\right|  ^{\alpha}\varepsilon^{-1/2}$ for any $\alpha>1/4$
would do.

Given an arbitrary subset $B\subset\Lambda,$ we write $N_{B}$ for the number
of big cells containing sites of $B$. If $A\subset\Lambda$ is another subset,
then we write $\mathcal{N}_{B}(A)$ for the number of those cells containing
sites of $B$ but no site of $A$ or of $\mathbb{Z}^{2}\setminus\Lambda$. We
shortly write $\mathcal{N}_{B}=\mathcal{N}_{B}(\mathcal{A}),$ when
$\mathcal{A}$ is our standard random subset, distributed according to
$\nu_{\Lambda}^{a,b}.$ Let $\rho=\left.  \left|  \log\varepsilon\right|
^{-2}\varepsilon\,\left|  B\right|  \right/  2N_{B}$. We want to prove that
\begin{equation}
\nu_{\Lambda}^{a,b}(\mathcal{N}_{B}>\rho\,N_{B})\leq\mathrm{e}^{-C\,\left|
\log\varepsilon\right|  \,^{-1/2}\;{\varepsilon}\;\left|  B\right|
}\,,\label{eq_manydirty}%
\end{equation}
provided $\varepsilon$ is small enough (independently of $B$). Notice that
\[
\tfrac{1}{2}\left|  \log\varepsilon\right|  ^{-2}{\varepsilon}\leq\rho
\leq\tfrac{1}{2}\,.
\]
(\ref{eq_manydirty}) is an easy consequence of Proposition~\ref{pro_fat}.
Indeed, we can apply the latter to get
\begin{align*}
\nu_{\Lambda}^{a,b}(\mathcal{N}_{B}>\rho\,N_{B}) &  \leq\sum_{k>\rho\,N_{B}%
}^{N_{B}}\binom{N_{B}}{k}\exp\left[  -C\,\left|  \log\varepsilon\right|
^{-1/2}{\varepsilon}\,k{\varepsilon}^{-1}\left|  \log\varepsilon\right|
^{2}\right] \\
&  =\sum_{k>\rho\,N_{B}}^{N_{B}}\binom{N_{B}}{k}\exp\left[  -C\,\left|
\log\varepsilon\right|  ^{3/2}\,k\right] \\
&  \leq\inf_{t\geq0}\left\{  \mathrm{e}^{-t\rho\,N_{B}}\sum_{k=0}^{N_{B}%
}\binom{N_{B}}{k}\exp\left[  (t-C\,\left|  \log\varepsilon\right|
^{3/2}k\right]  \right\} \\
&  =\inf_{t\geq0}\left\{  \mathrm{e}^{-t\rho}\left[  1+\exp\left[
t-C\,\left|  \log\varepsilon\right|  ^{3/2}\right]  \right]  \right\}
^{N_{B}}\\
&  \leq\exp\{-\tfrac{1}{2}C\,\left|  \log\varepsilon\right|  ^{3/2}%
\,\rho\,N_{B}\}\\
&  =\exp\left[  -\tfrac{1}{4}C\,\left|  \log\varepsilon\right|  ^{-1/2}%
\,{\varepsilon}\,\left|  B\right|  \right]  \,.
\end{align*}

\subsubsection{Step 3: Arbitrary sets}

Let now $B$ be an arbitrary subset of $\Lambda$. By~\eqref{eq_manydirty}, we
know that
\[
\nu_{\Lambda}^{a,b}(\mathcal{A}\cap B=\emptyset)\leq\nu_{\Lambda}^{a,b}\left(
\mathcal{A}\cap B=\emptyset\,\left|  \,\mathcal{N}_{B}<\frac{1}{2}\left|
\log\varepsilon\right|  ^{-2}\,\varepsilon\,\left|  B\right|  \right.
\right)  +\exp\left[  -C\,\left|  \log\varepsilon\right|  ^{-1/2}%
\,\varepsilon\,\left|  B\right|  \right]  .
\]
In order to finish the proof of the theorem, it remains to estimate the first
summand on the right-hand side. The idea is to essentially repeat the argument
used in the proof of Proposition~\ref{pro_fat}, using the fact that there are
already many pinned sites close to $B$. Let us therefore suppose, without loss
of generality, that
\[
\left\{  A:A\cap B=\emptyset,\,\mathcal{N}_{B}(A)<\frac{1}{2}\left|
\log\varepsilon\right|  ^{-2}\varepsilon\left|  B\right|  \right\}
\neq\emptyset
\]
(otherwise the conditional probability is simply $0$ and there is nothing to
prove). Then we have, as in (\ref{eq_probclean}),
\begin{align*}
&  \nu_{\Lambda}^{a,b}\left(  \mathcal{A}\cap B=\emptyset\left|
\,\mathcal{N}_{B}^{{}}<\frac{1}{2}\left|  \log\varepsilon\right|
^{-2}\,\varepsilon\,\left|  B\right|  \right.  \right) \\
&  \leq\left\{  \inf_{A}\sum_{D\subset B}(\mathrm{e}^{b}-1)^{\left|  D\right|
}\mu_{\Lambda}\left(  \left|  \phi_{x}\right|  \leq a,\,\forall x\in
D\bigm\vert\,\left|  \phi_{x}\right|  \leq a,\,\forall x\in A\right)
\right\}  ^{-1}\\
&  \leq\left\{  \inf_{A}\sum_{D\subset B^{\mathrm{g}}(A)}(\mathrm{e}%
^{b}-1)^{\left|  D\right|  }\mu_{\Lambda}\left(  \left|  \phi_{x}\right|  \leq
a,\,\forall x\in D\bigm\vert\,\left|  \phi_{x}\right|  \leq a,\,\forall x\in
A\right)  \right\}  ^{-1}%
\end{align*}
where the infimum is taken over sets $A$ with $A\cap B=\emptyset$ and
$\mathcal{N}_{B}(A)<\frac{1}{2}\left|  \log\varepsilon\right|  ^{-2}%
\varepsilon\left|  B\right|  ,$ and where $B^{\mathrm{g}}(A)$ is the set of
``good'' points in $B:$ those sharing a big box with at least one point from
$A$ or $\mathbb{Z}^2\setminus \Lambda$. It is easy to estimate the inner
probability. Indeed, numbering the elements of $D=\{t_{1},\ldots,t_{\left| 
D\right|  }\}$, we can write
\begin{align*}
\mu_{\Lambda}\left(  \left.  \left|  \phi_{x}\right|  \leq a,\forall x\in
D\right|  \,\left|  \phi_{x}\right|  \leq a,\,\forall x\in A\right)   &
=\prod_{k=1}^{\left|  D\right|  -1}\mu_{\Lambda}\left(  \left.  \left|
\phi_{t_{k+1}}\right|  \leq a\right|  \,\left|  \phi_{x}\right|  \leq
a,\,\forall x\in A\cup\left\{  t_{1},\ldots,t_{k}\right\}  \right) \\
&  \geq\prod_{k=1}^{\left|  D\right|  -1}\frac{1}{2}\left(  \frac{a}%
{4\mu_{A^{c}\backslash\left\{  t_{1},\ldots,t_{k}\right\}  }\left(  \left|
\phi_{t_{k+1}}\right|  \right)  }\wedge\frac{1}{2}\right)
\end{align*}
where the last inequality follows from Lemmas~5.4 and~5.5 of~\cite{DV}. The
expected value is easily estimated using the random walk representation:
\begin{align*}
\mu_{A^{c}\backslash\left\{  t_{1},\ldots,t_{k}\right\}  }\left(  \left|
\phi_{t_{k+1}}\right|  \right)   &  \leq\sqrt{\mu_{A^{c}\backslash\left\{
t_{1},\ldots,t_{k}\right\}  }\left(  \phi_{t_{k+1}}^{2}\right)  }\\
&  \leq\sqrt{\mu_{A^{c}\backslash\left\{  t_{1},\ldots,t_{k}\right\}  }^{\ast
}\left(  \phi_{t_{k+1}}^{2}\right)  }\\
&  \leq C\sqrt{\left|  \log\varepsilon\right|  }%
\end{align*}
where the second inequality follows from Brascamp-Lieb, and the last one
follows from~\eqref{eq_GF_one_obstacle}, since the last probability is bounded
by the Green function of the random walk killed as it hits the closest site of
$A$ or of $\mathbb{Z}^{2}\setminus\Lambda$ located in the same cell as
$t_{k+1}$ (there is such a site since $t_{k+1}\in B^{\mathrm{g}}(A)$).

Therefore,
\[
\mu_{\Lambda}\left(  \left.  \left|  \phi_{x}\right|  \leq a,\,\forall x\in
D\right|  \,\left|  \phi_{x}\right|  \leq a,\,\forall x\in A\right)
\geq\left(  \frac{C\left(  2a\wedge1\right)  }{\sqrt{\left|  \log
\,\varepsilon\right|  }}\right)  ^{\left|  D\right|  }.
\]
This finally yields
\begin{align*}
\nu_{\Lambda}^{a,b}(\mathcal{A}\cap B=\emptyset\,|\,\mathcal{N}_{B}^{{}%
}<\tfrac{1}{2}\left|  \log\varepsilon\right|  ^{-2}\varepsilon\,\left|
B\right|  ) &  \leq\left\{  \inf_{A}\sum_{D\subset B^{\mathrm{g}}}\left(
\frac{D\,\varepsilon}{\sqrt{\left|  \log\varepsilon\right|  }}\right)
^{\left|  D\right|  }\right\}  ^{-1}\\
&  \leq\exp\left[  -C\,\left|  \log\varepsilon\right|  ^{-1/2}\,{\varepsilon
}\,\inf_{A}\left|  B^{g}(A)\right|  \right]  \,.
\end{align*}
The conclusion follows easily since
\[
\left|  B^{g}(A)\right|  \geq\left|  B\right|  -\mathcal{N}_{B}(A)\left|
\log\varepsilon\right|  ^{2}\,\varepsilon^{-1}\geq\left|  B\right|  /2,
\]
when $\mathcal{N}_{B}(A)<\frac{1}{2}\left|  \log\varepsilon\right|
^{-2}\,\varepsilon\,\left|  B\right|  .$

\subsection{Proof of point 4.}

We assume here $d\geq3.$ The desired inequality follows from \eqref{eq_rem}
and, using Lemmas~5.4 and 5.5 of~\cite{DV},
\begin{align*}
\frac{Z_{\Lambda}(A\cup\{x\})}{Z_{\Lambda}(A)}=\mu_{\Lambda}(|\phi_{x}|\leq
a\,|\,|\phi_{y}|\leq a\;\forall y\in A) &  \geq\tfrac{1}{2}(\frac{a}%
{4\mu_{\Lambda}(|\phi_{x}|)}\wedge\tfrac{1}{2})\\
&  \geq\tfrac{1}{2}(\frac{a}{4\sqrt{\mu(|\phi_{x}|^{2})}}\wedge\tfrac{1}%
{2})\geq C_{3}^{\prime}(2a\wedge1)\,.
\end{align*}

\section{Asymptotics of the variance \label{SectVariance}}

\subsection{Proof of Theorem \ref{thm_lowerbound}, and the lower bound in
Theorem \ref{thm_critical}}

Let $\Lambda$ be a square in $\mathbb{Z}^{2}$, centered at the origin, and
with large enough sidelength (the thermodynamic limit is taken at the end).
Let $B_{\overline{\mathfrak{e}}}(0)\overset{\mathrm{def}}{=}\left\{
x\in\mathbb{Z}^{2}:\left\|  x\right\|  _{\infty}\leq\tfrac12 {\overline
{\mathfrak{e}}}^{-1/2} \left\vert \log\overline{\mathfrak{e}} \right\vert
^{-1/4}\right\}  $. Using \eqref{eq_exp_sw}, we get
\begin{align*}
\mu_{\Lambda}^{a,b}(\phi_{0}^{2})  &  =\sum_{A\subset\Lambda}\nu_{\Lambda
}^{a,b}(A)\;\mu_{\Lambda}(\phi_{0}^{2}\,|\,\left|  \phi_{x}\right|  \leq
a,\,\forall x\in A)\\
&  \geq\nu_{\Lambda}^{a,b}(\mathcal{A}\cap B_{\overline{\mathfrak{e}}
}(0)=\emptyset)\;\inf_{A\cap B_{\overline{\mathfrak{e}}}(0)=\emptyset}%
\mu_{\Lambda}(\phi_{0}^{2}\,|\,\left|  \phi_{x}\right|  \leq a,\,\forall x\in
A)\,.
\end{align*}
Using $\mu_{\Lambda}(\phi_{0}^{2}\,|\,\left|  \phi_{x}\right|  \leq
a,\,\forall x\in A)\geq\tfrac{1}{2}\mu_{A^{c}}(\phi_{0}^{2})$ (see~ \cite{DV},
Lemma~5.4), we get
\[
\mu_{\Lambda}^{a,b}(\phi_{0}^{2})\geq\tfrac{1}{2}\nu_{\Lambda}^{a,b}
(\mathcal{A}\cap B_{\overline{\mathfrak{e}}}(0)=\emptyset)\;\inf_{A\cap
B_{\overline{\mathfrak{e}}}(0)=\emptyset}\mu_{A^{c}}(\phi_{0}^{2})\,.
\]
Now, by part 1 of Theorem~\ref{thm_obstacles}, we have $\nu_{\Lambda}
^{a,b}(\mathcal{A}\cap B_{\overline{\mathfrak{e}}}(0)=\emptyset)\geq
1-2C\,\left|  \log\underline{\mathfrak{e}}\right|  ^{-1}$, and by the inverse
Brascamp-Lieb inequality~\cite{DGI}, $\mu_{A^{c}}(\phi_{0}^{2})\geq\tfrac
{1}{\overline{c} }\mu_{A^{c}}^{\star}(\phi_{0}^{2})=\tfrac{1}{\beta
\overline{c}}G_{A}(0,0)$, where the last quantity is the Green function for
the simple random walk killed as it enters the set $A$. Clearly $G_{A}(0,0)$
is minimum when $A=\mathbb{Z}^{2}\setminus B_{\overline{\mathfrak{e}}}(0)$.
The conclusion follows from \eqref{eq_G_in_B}.

The case of $\delta$-pinning is identical. Note however that in this case we
do not lose the factor $\tfrac{1}{2}$ since the pinned sites are already at
$0$. Also, in the Gaussian case we do not need Brascamp-Lieb, and therefore we
do not get the factor $\tfrac{1}{\overline{c}}$ either. Therefore, we get in
this case
\[
\mu^{\ast,\varepsilon}\left(  \phi_{0}^{2}\right)  \geq\frac{\left|
\log\left(  \sqrt{\beta}\varepsilon\right)  \right|  }{2\pi\beta\sqrt
{\det\mathcal{Q}}}-C\log\left|  \log\left(  \sqrt{\beta}\varepsilon\right)
\right|  ,
\]
which proves the lower bound in Theorem \ref{thm_critical}.

\subsection{Proof of the upper bound in Theorem \ref{thm_critical}}

We apply Remark \ref{RemRescale}, and therefore assume $\beta=1.$ Using
Corollary \ref{CorDom}, we have
\begin{align}
\,\mu_{\Lambda}^{\star,\mathcal{\varepsilon}}(\phi_{0}^{2})  &  \leq
\sum_{n\geq0}\rho_{+}\otimes\mathrm{P}_{0}[X_{n}=0,\tau_{\mathcal{A}^{c}
}>n]\nonumber\\
&  =\sum_{n=0}^{n_{0}}\mathrm{P}_{0}[X_{n}=0]+\sum_{n>n_{0}}\rho_{+}
\otimes\mathrm{P}_{0}[\tau_{\mathcal{A}^{c}}>n]\nonumber\\
&  =G^{n_{0}}(0,0)+\sum_{n>n_{0}}\mathrm{E}_{0}[(1-p_{+})^{\left|
X_{[0,n]}\right|  }]\,\label{eq_22}%
\end{align}
where we choose $n_{0}=n_{0}(\mathcal{\varepsilon})=\mathcal{\varepsilon}%
^{-1}\,\left|  \log\varepsilon\right|  ^{\eta}$, for some $\eta>0$ to be chosen
later. Then the $n_{0}$-step Green
function in the right-hand side of the last equation has the following
asymptotics, see~\eqref{eq_G_upto_n},
\[
G^{n_{0}}(0,0)=(2\pi\sqrt{\det Q})^{-1}\left|  \log\varepsilon\right|
+O(\log\left| \log\varepsilon\right| )\,.
\]
The claim will be proved if we show that the second term in~\eqref{eq_22} does
not contribute more than $O(\log\left|  \log\varepsilon\right|  )$; we are
actually going to check that it is even $o(1)$ as $\mathcal{\varepsilon}$ goes
to zero. Indeed, introducing a small constant $\kappa>0$, it can be estimated
in the following way:
\[
\sum_{n>n_{0}}\mathrm{E}_{0}[(1-p_{+})^{\left|  X_{[0,n]}\right|  }]\leq
\sum_{n>n_{0}}(1-p_{+})^{\kappa\,n/\log n}+\sum_{n>n_{0}}\mathrm{P}
_{0}[\left|  X_{[0,n]}\right|  \leq\kappa\,n/\log n]\,.
\]
By Proposition \ref{PropModerate}, we see that $\mathrm{P}_{0}[\left|
X_{[0,n]}\right|  \leq\kappa\,n/\log n]\leq n^{-2}$ provided $\kappa$ is
chosen small enough; this shows that the last sum is $o(1)$. To see that this
is also true for the first one, we bound it as follows (remember that
$n_{0}\rightarrow\infty$ when $\mathcal{\varepsilon}\rightarrow0$):
\begin{align*}
\sum_{n>n_{0}}(1-p_{+})^{\kappa\,n/\log n}  &  \leq\sum_{n>n_{0}}
e^{-p_{+}\kappa\,n/\log n}\\
&  \leq\int_{n_{0}-1}^{\infty}e^{-p_{+}\kappa\,x/\log x}\,\mathrm{d}x\\
&  \leq\int_{n_{0}/2}^{\infty}e^{-\tfrac{1}{2}p_{+}\kappa\,y}\,\mathrm{d}y\\
&  =\frac{2}{p_{+}\kappa}\,e^{-\tfrac{1}{4}p_{+}\kappa\,n_{0}}%
\end{align*}
which is $o(1)$ by definition of $p_{+}$ and $n_{0}$, provided we take $\eta$
sufficiently large (depending on $\kappa$).

\section{Asymptotics of the mass: Proof of Theorem~\ref{thm_critical2}%
\label{SectMass}}

We discuss the $2$-dimensional case in details. The simpler higher-dimensional
case follows exactly in the same way by using Theorem \ref{thm_obstacles}
parts 1 and 4 instead of parts 2 and 3.

We consider $x\in\mathbb{Z}^{2}$ sufficiently far away from $0$. We take
$\Lambda$ to be a finite box in $\mathbb{Z}^{2}$, and prove the estimates when
$\Lambda$ is large enough, depending possibly on $x.$ This then proves the
estimates in the thermodynamic limit. Remember that we assume here that
$\left(  p\left(  x\right)  \right)  _{x\in\mathbb{Z}^{2}}$ has an exponential
moment. Furthermore, we assume that $p$ is irreducible and aperiodic.

\subsection*{Proof of the upper bound}

We denote by $\mathrm{E}_{x,y}^{(n)}$ the expectation for the random walk
starting in $x$ and conditioned on $X_{n}=y$, provided the probability of the
latter event is positive. Using the Corollary \ref{CorDom}, we have
\begin{align*}
\mu^{\ast,\varepsilon}\left(  \phi_{0}\phi_{x}\right)   &  \geq\sum_{n\geq
0}\mathrm{E}_{0}\left(  \exp\left[  -C\varepsilon|\log\varepsilon
|^{-1/2}\left|  X_{[0,n]}\right|  \right]  ;\,X_{n}=x\right) \\
&  =\sum_{n=0}^{\infty}p_{n}(x)\mathrm{E}_{0}\left(  \left.  \exp\left[
-C\varepsilon|\log\varepsilon|^{-1/2}\left|  X_{[0,n]}\right|  \right]
\right|  X_{n}=x\right) \\
&  \geq\sum_{n=0}^{\infty}p_{n}(x)\exp\left[  -C\varepsilon|\log
\varepsilon|^{-1/2}\mathrm{E}_{0,x}^{(n)}\left(  \left|  X_{[0,n]}\right|
\right)  \right] \\
&  \geq p_{m}(x)\exp\left[  -C\varepsilon|\log\varepsilon|^{-1/2}%
\mathrm{E}_{0,x}^{(m)}\left(  \left|  X_{[0,m]}\right|  \right)  \right]  ,
\end{align*}
where
\[
m=m(\left|  x\right|  ,\varepsilon)\overset{\mathrm{def}}{=}\left[  \left|
\log\varepsilon\right|  ^{3/4}\varepsilon^{-1/2}\left|  x\right|  \right]  .
\]
We apply Proposition \ref{PropTieddown}, and use
\[
p_{m}(x)\geq\frac{C}{m}\exp\left[  -mI\left(  \frac{x}{m}\right)  \right]
\geq\frac{C}{m}\exp\left[  -\rho\frac{\left|  x\right|  ^{2}}{m}\right]  ,
\]
for some positive $\rho$, see~Proposition~\ref{PropSaddle}. So we get
\begin{align*}
&  \sum_{n=0}^{\infty}\mathrm{E}_{0}\left(  \exp\left[  -C\varepsilon
|\log\varepsilon|^{-1/2}\left|  X_{[0,n]}\right|  \right]  ;X_{n}=x\right) \\
&  \geq\frac{C^{\prime\prime}}{\left|  \log\varepsilon\right|  ^{3/4}%
\varepsilon^{-1/2}\left|  x\right|  }\exp\left[  -\rho\left|  \log
\varepsilon\right|  ^{-3/4}\varepsilon^{1/2}\left|  x\right|  -C^{\prime
}\varepsilon|\log\varepsilon|^{-1/2}\frac{\left|  \log\varepsilon\right|
^{3/4}\varepsilon^{-1/2}\left|  x\right|  }{\log(\left|  \log\varepsilon
\right|  ^{3/4}\varepsilon^{-1/2})}\right] \\
&  \geq\exp\left[  -C^{\prime\prime\prime}\left|  \log\varepsilon\right|
^{-3/4}\varepsilon^{1/2}\left|  x\right|  \right]
\end{align*}
for small enough $\varepsilon>0,$ and then large enough $\left|  x\right|  $.
This proves the lower bound.

There is a trivial modification necessary for $d\geq3:$ We have to replace the
use of Proposition \ref{PropTieddown} by the completely trivial bound $\bigl
|X_{[0,n]}\bigr|
\leq n+1.$

\subsection*{Proof of the lower bound}

We start by proving that the logarithmic asymptotics for the 2-point function
$\mu_{\Lambda}^{\ast,\varepsilon}(\phi_{0}\phi_{x})$ are entirely determined
by the probability that the random walk reaches $x$ before dying.

\begin{lemma}
1.
\[
\nu^{\varepsilon}\otimes\mathrm{P}_{0}[\mathrm{T}_{\{x\}}<\mathrm{T}%
_{\mathcal{A}}]\overset{\mathrm{def}}{=}\lim_{\Lambda\nearrow\mathbb{Z}^{2}%
}\nu_{\Lambda}^{\varepsilon}\otimes\mathrm{P}_{0}[\mathrm{T}_{\{x\}}%
<\tau_{\mathcal{A}^{c}}]
\]
exists for all $x\in\mathbb{Z}^{d}$.\newline 2. For all $x\in\mathbb{S}^{1}$,
\begin{align*}
\limsup_{k\rightarrow\infty}\frac{1}{k}\log\mu^{\ast,\varepsilon}(\phi_{0}%
\phi_{\lbrack kx]}) &  \leq\limsup_{k\rightarrow\infty}\frac{1}{k}\log
\nu^{\varepsilon}\otimes\mathrm{P}_{0}[\mathrm{T}_{\{[kx]\}}<\tau
_{\mathcal{A}^{c}}]\,,\\
\liminf_{k\rightarrow\infty}\frac{1}{k}\log\mu^{\ast,\varepsilon}(\phi_{0}%
\phi_{\lbrack kx]}) &  \geq\liminf_{k\rightarrow\infty}\frac{1}{k}\log
\nu^{\varepsilon}\otimes\mathrm{P}_{0}[\mathrm{T}_{\{[kx]\}}<\tau
_{\mathcal{A}^{c}}]\,.
\end{align*}
(That these limits actually exist is proved in Appendix~\ref{app_mass}.) \label{lem_sameasympt}
\end{lemma}

\begin{proof}
\noindent1. If $\Lambda^{\prime}\subset\Lambda\Subset\mathbb{Z}^{2}$, FKG
property of $\nu_{\Lambda}^{\mathcal{\varepsilon}}$ implies $\nu_{\Lambda
}^{\mathcal{\varepsilon}}(\mathcal{A}\cap D=\emptyset)\geq\nu_{\Lambda
^{\prime}}^{\mathcal{\varepsilon}}(\mathcal{A}\cap D=\emptyset)$ for any set
$D$, see~\eqref{eq_monotone}. Therefore
\begin{align*}
\nu_{\Lambda}^{\mathcal{\varepsilon}}\otimes\mathrm{P}_{0}[\mathrm{T}%
_{\{x\}}<\tau_{\Lambda\setminus\mathcal{A}}] &  =\mathrm{E}_{0}[\mathrm{I}%
\left(  \mathrm{T}_{\left\{  x\right\}  }<\tau_{\Lambda}\right)
\,\nu_{\Lambda}^{\mathcal{\varepsilon}}(\mathcal{A}\cap X_{[0,\mathrm{T}%
_{\{x\}}]}=\emptyset)]\\
&  \geq\mathrm{E}_{0}[\mathrm{I}\left(  \mathrm{T}_{\left\{  x\right\}  }%
<\tau_{\Lambda^{\prime}}\right)  \,\nu_{\Lambda^{\prime}}%
^{\mathcal{\varepsilon}}(\mathcal{A}\cap X_{[0,\mathrm{T}_{\{x\}}]}%
=\emptyset)]\\
&  =\nu_{\Lambda^{\prime}}^{\mathcal{\varepsilon}}\otimes\mathrm{P}%
_{0}[\mathrm{T}_{\{x\}}<\tau_{\Lambda^{\prime}\setminus\mathcal{A}}]\,,
\end{align*}
which proves the claim since the probabilities are bounded by $1$.

\medskip\noindent2. Using the expansion~\eqref{eq_exp_pinning}, we can write
\begin{align*}
\,\mu_{\Lambda}^{\star,\varepsilon}(\phi_{0}\,\phi_{x}) &  =\sum_{n\geq0}%
\nu_{\Lambda}^{\varepsilon}\otimes\mathrm{P}_{0}[X_{n}=x,X_{[0,n]}%
\subset\mathcal{A}^{c}]\\
&  =\sum_{n\geq0}\sum_{A\subset\Lambda}\nu_{\Lambda}^{\varepsilon
}(A)\;\mathrm{P}_{0}[\mathrm{T}_{\{x\}}<\tau_{A^{c}}]\;\mathrm{P}_{x}%
[X_{n}=x,X_{[0,n]}\subset A^{c}]\\
&  =\sum_{A\subset\Lambda}\nu_{\Lambda}^{\varepsilon}(A)\;\mathrm{P}%
_{0}[\mathrm{T}_{\{x\}}<\tau_{A^{c}}]\;G_{A}(x,x)\\
&  =\sum_{R\geq0}\sum_{\substack{A\subset\Lambda\\ d_{\infty}(A,x)=R}%
}\nu_{\Lambda}^{\varepsilon}(A)\;\mathrm{P}_{0}[\mathrm{T}_{\{x\}}<\tau
_{A^{c}}]\;G_{A}(x,x)\\
&  \leq\nu_{\Lambda}^{\varepsilon}\otimes\mathrm{P}_{0}[\mathrm{T}%
_{\{x\}}<\tau_{\mathcal{A}^{c}}]\;\sum_{R\geq0}\max_{y:\Vert x-y\Vert
_{\scriptscriptstyle\infty}=R}G_{\mathbb{Z}^{2}\setminus\{y\}}(x,x)\\
&  \hspace*{5cm}\times\nu_{\Lambda}^{\varepsilon}\otimes\mathrm{P}%
_{0}[d_{\infty}(\mathcal{A},x)=R\,|\,\mathrm{T}_{\{x\}}<\tau_{\mathcal{A}^{c}%
}]\\
&  \leq\nu_{\Lambda}^{\varepsilon}\otimes\mathrm{P}_{0}[\mathrm{T}%
_{\{x\}}<\tau_{\mathcal{A}^{c}}]\;\sum_{R\geq0}C\log R\;\nu_{\Lambda
}^{\varepsilon}\otimes\mathrm{P}_{0}[d_{\infty}(\mathcal{A}%
,x)=R\,|\,\mathrm{T}_{\{x\}}<\tau_{\mathcal{A}^{c}}]\,.
\end{align*}
We therefore have to bound the conditional probability. This can be done as
follows:
\[
\nu_{\Lambda}^{\varepsilon}\otimes\mathrm{P}_{0}[d_{\infty}(\mathcal{A}%
,x)=R\,|\,\mathrm{T}_{\{x\}}<\tau_{\mathcal{A}^{c}}]\leq\frac{\nu_{\Lambda
}^{\varepsilon}\otimes\mathrm{P}_{0}[d_{\infty}(\mathcal{A},x)=R]}%
{\nu_{\Lambda}^{\varepsilon}\otimes\mathrm{P}_{0}[\mathrm{T}_{\{x\}}%
<\tau_{\mathcal{A}^{c}}]}\wedge1\leq\frac{e^{-C(\varepsilon)R^{2}}%
}{e^{-C^{\prime}(\varepsilon)|x|^{2}}}\wedge1\,,
\]
where we used Theorem~\ref{thm_obstacles} to bound the numerator and the bound
on the denominator follows from
\[
\nu_{\Lambda}^{\varepsilon}\otimes\mathrm{P}_{0}[\mathrm{T}_{\{x\}}%
<\tau_{\mathcal{A}^{c}}]\geq\sum_{n\geq0}\mathrm{E}_{0}[\,\mathrm{I}%
(\mathrm{T}_{\{x\}}=n)\,(1-p_{-})^{|X_{[0,n]}|}]\geq(1-p_{-})^{|x|^{2}%
}\mathrm{P}_{0}[\mathrm{T}_{\{x\}}\leq x^{2}-1]\,,
\]
and the local CLT. Therefore the sum over $R$ is smaller than $C(\varepsilon
)(|x|\log|x|+1)$, which proves the first claim.

\medskip To prove the second claim, notice that
\[
\,\mu^{\star,\varepsilon}_{\Lambda}(\phi_{0}\,\phi_{x}) = \sum_{R\geq0}
\sum_{\substack{A\subset\Lambda\\ d_{\infty}(A,x)=R}} \nu^{\varepsilon
}_{\Lambda}(A)\; \mathrm{P}_{0}[\mathrm{T}_{\{x\}} < \tau_{A^{c}}]\;
G_{A}(x,x) \geq\nu^{\varepsilon}_{\Lambda}\otimes\mathrm{P}_{0} [\mathrm{T}%
_{\{x\}} < \tau_{\mathcal{A}^{c}}]\,,
\]
since $G_{A}(x,x)\geq1$ (one can restrict the sum over sets $A$ not containing
$x$, since otherwise the probability of reaching $x$ is 0).
\end{proof}

Let $\Delta(\varepsilon)=\varepsilon^{-1/2} |\log\varepsilon|^{3/4}$. We
consider a partition of $\mathbb{Z}^{2}$ into cells of width $\Delta
(\varepsilon)$, and write, for $y\in\mathbb{Z}^{2}$, $B_{y}$ for the cell
containing $y$, and $\overline{B}_{y}$ for the square composed of
$(2M+1)\times(2M+1)$ cells with middle-cell $B_{y}$, where $M$ is a big
integer to be chosen later. We introduce the following stopping-times

\begin{itemize}
\item $T_{0} = 0$ ;

\item $T_{k} = \min\{ n>T_{k-1}\,:\, \overline{B}_{X_{n}} \cap\overline
{B}_{X_{T_{l}}} = \emptyset\quad\forall l<k\} \quad\quad(k\geq1)$ ;

\item $T_{k}^{\prime}=\min\{n>T_{k}\,:\,X_{n}\not \in B_{M\Delta(\varepsilon
)}(X_{T_{k}})\}\quad\quad(k\geq0)$,
\end{itemize}

where $B_{M\Delta(\varepsilon)}(y)$ is the ball of radius $M\Delta
(\varepsilon)$ and center $y$. Let also $\overline{k}=\max\{k\geq
0\,:\,T_{k}<\mathrm{T}_{\overline{B}_{x}}\}$, and let $c$ be some small
constant to be chosen later. We then have (remember that $p_{+}=C\varepsilon
|\log\varepsilon|^{-1/2}$)
\begin{align*}
\nu^\varepsilon\otimes\mathrm{P}_{0}[\mathrm{T}_{\{x\}}<T_{\mathcal{A}}] &
\leq \mathrm{E}_{0}\bigl[(1-p_{+})^{|X_{[0,\mathrm{T}_{x}]}|}\bigr]\\
&  \leq\mathrm{E}_{0}\bigl[(1-p_{+})^{|X_{[0,\mathrm{T}_{x}]}|}\;\mathrm{I}%
(\overline{k}> c|x|/\Delta(\varepsilon))\bigr]+\mathrm{P}_{0}\bigl
[\overline{k}\leq c|x|/\Delta(\varepsilon)\bigr]%
\end{align*}
By Proposition~\ref{PropCross}, the last probability is bounded from above by
$\mathrm{e}^{-C|x|}$, with $C>0$ independent of $\varepsilon$, provided $c$ is
chosen small enough; indeed the total number of cells visited is certainly
smaller than $(2M+1)^{2}(\overline{k}+1)$. Let us now consider the first
term. Clearly,
\begin{align*}
\mathrm{E}_{0}\bigl[(1-p_{+})^{|X_{[0,\mathrm{T}_{x}]}|}\;\mathrm{I}%
(\overline{k}> c|x|/\Delta(\varepsilon))\bigr] &  \leq\mathrm{E}_{0}\bigl
[\prod_{k=0}^{\overline{k}-1}(1-p_{+})^{|X_{[T_{k},T_{k}^{\prime}]}%
|}\;\mathrm{I}(\overline{k}\geq c|x|/\Delta(\varepsilon))\bigr]\\
&  \leq\mathrm{E}_{0}\bigl[\prod_{k=0}^{c|x|/\Delta(\varepsilon)}%
(1-p_{+})^{|X_{[T_{k},T_{k}^{\prime}]}|}\bigr]\\
&  \leq\left\{  \mathrm{E}_{0}\bigl[(1-p_{+})^{|X_{[0,\tau_{B_{M\Delta
(\varepsilon)}(0)}]}|}\bigr]\right\}  ^{c|x|/\Delta(\varepsilon)}\,,
\end{align*}
The conclusion follows, since the latter expectation can easily be bounded.
Choosing some $C_{1}>0,$ we split as follows:
\begin{align*}
\mathrm{E}_{0}\bigl[(1-p_{+})^{|X_{[0,\tau_{B_{M\Delta(\varepsilon)}(0)}]}%
|}\bigr]&\leq e^{-p_{+}C_{1}/p_{+}}+\mathrm{P}_{0}[\tau_{B_{M\Delta
(\varepsilon)}(0)}<M\Delta(\varepsilon)^{2}]\\
&\hspace{5cm}+\mathrm{P}_{0}[|X_{[0,M\Delta
(\varepsilon)^{2}]}|\leq C_{1}/p_{+}]\\
&\leq 3/4\,,
\end{align*}
We now choose first $C_{1}$ such that the first summand is $\leq1/4.$ Next,
observe that by the invariance principle for the random walk, we have
\[
\lim_{M\rightarrow\infty}\mathrm{P}_{0}[\tau_{B_{M\Delta(\varepsilon)}%
(0)}<M^{2}\Delta(\varepsilon)^{2}]=\mathrm{P}_{0}\left[  \sigma_{B_{1}\left(
0\right)  }<1\right]  ,
\]
where $\sigma$ is the exit time of a Brownian motion with covariance
$\mathcal{Q},$ and therefore
$$
\lim_{M\rightarrow\infty}\mathrm{P}_{0}
[\tau_{B_{M\Delta(\varepsilon)}(0)}<M\Delta(\varepsilon)^{2}]=0\,,
$$
uniformly in $\varepsilon\leq1/2,$ say. As $M\rightarrow\infty,$ also the third
summand is
converging to $0,$ uniformly in $\varepsilon\leq1/2,$ which follows from the
law of large numbers for the range of the random-walk, see~ \cite{JP} (notice
that $C_{2}/p_{+}=C_{2}\,C\,\Delta_{d}(\varepsilon)^{2}/\log\Delta
_{d}(\varepsilon)$). Therefore, we can make the second and the third summand
$\leq1/4$ by choosing $M$ appropriately.

Again, there is a trivial modification for $d\geq3:$ One chooses
$\Delta(\varepsilon)=\varepsilon^{-1/2}$ and appeal to Theorem T1.4.1 of
~\cite{Spitzer} for the law of large number for the range of these transient
random walks; remember that $p_{+}$ is now equal to $C\varepsilon$.

\appendix

\section{Existence of the mass}

\label{app_mass}

In this appendix, we prove existence of the limit in~\eqref{eq_mass}, in the
case of Gaussian interactions and $\delta$-pinning. According to
Lemma~\ref{lem_sameasympt}, the existence of the mass is a consequence of the
following result.

\begin{lemma}%
\[
\lim_{k\rightarrow\infty}-\frac{1}{k\left|  x\right|  }\,\log\nu^{\varepsilon
}\otimes\mathrm{P}_{0} [\mathrm{T}_{\{x\}}<\mathrm{T}_{\mathcal{A}}] =
\sup_{k}\,-\frac{1}{k\left|  x\right|  }\,\log\nu^{\varepsilon}\otimes
\mathrm{P}_{0} [\mathrm{T}_{\{x\}}<\mathrm{T}_{\mathcal{A}}]
\]
exists for all $x\in\mathbb{Z}^{2}$.
\end{lemma}

\begin{proof}
This follows from a standard subadditivity argument, since
\begin{align*}
\nu^{\mathcal{\varepsilon}} \otimes\mathrm{P}_{0} [\mathrm{T}_{\{[(k+l)x]\}}%
<\mathrm{T}_{\mathcal{A} }]  &  \geq\nu^{\mathcal{\varepsilon}}\otimes
\mathrm{P}_{0}[\mathrm{T}_{\{[kx]\}} <\mathrm{T}_{\{[(k+l)x]\}}^{\prime
}<\mathrm{T}_{\mathcal{A}}]\\
&  =\mathrm{E}_{0}[\mathrm{I}\left(  \mathrm{T}_{\left\{  [kx]\right\}
}<\mathrm{T}_{\left\{  [( k+l) x]\right\}  }^{\prime}<\infty\right)
\,\nu^{\mathcal{\varepsilon}}(\mathcal{A}\cap X_{[0,\mathrm{T}_{\{[(k+l)x]\}}%
^{\prime} ]}=\emptyset)]\\
&  \geq\mathrm{E}_{0}[\mathrm{I}\left(  \mathrm{T}_{\left\{  [kx]\right\}
}<\mathrm{T}_{\left\{  [( k+l) x]\right\}  }^{\prime}<\infty\right)
\,\nu^{\mathcal{\varepsilon}}(\mathcal{A}\cap X_{[0,\mathrm{T}_{\{[kx]\}}%
]}=\emptyset)\\
&  \hspace*{5cm}\times\nu^{\mathcal{\varepsilon}}(\mathcal{A}\cap
X_{[\mathrm{T}_{\{[kx]\}},\mathrm{T}_{\{[(k+l)x]\}}^{\prime}]}=\emptyset)]\\
&  = C \nu^{\mathcal{\varepsilon}} \otimes\mathrm{P}_{0} [\mathrm{T}%
_{\{[kx]\}}<\mathrm{T}_{\mathcal{A} }]\;\nu^{\mathcal{\varepsilon}}
\otimes\mathrm{P}_{0} [\mathrm{T}_{\{[lx]\}}<\mathrm{T}_{\mathcal{A} }]\,,
\end{align*}
where $\mathrm{T}_{\left\{  [( k+l) x]\right\}  }^{\prime}\overset
{\mathrm{def}}{=}\min\{n>\mathrm{T}_{\{[kx]\}} \,|\, X_{n}=[(k+l)x]\}$, and
the inequality follows from the FKG property. The constant $C$, which depends
only on $p$, takes care of the possible discrepancy between $[(k+l)x]$ and
$[kx]+[lx]$.
\end{proof}

\section{Some properties of random walks\label{SectRW}}

\label{app_RW}

We keep the assumptions on $p$ made in the introduction. Especially, we always
assume the existence of a moment of order $2+\delta$ (\ref{CondH1}) and that
the random walk is irreducible and aperiodic.

We always use $C,$ $C^{\prime}$ for positive constants, not necessarily the
same at different occurrences, which may depend on $p(\,\cdot\,) and $d$,$ but
on nothing else.

\subsection{Properties of Green functions for random walks in dimension 2}

We denote by $a(x)=\sum_{n\geq0}(\mathrm{P}_{0}(X_{n}=0)-\mathrm{P}_{0}
(X_{n}=x))$ the potential kernel associated to the random walk.

For any $B\subset\mathbb{Z}^{2}$, $G_{B}(x,y)\overset{\mathrm{def}}%
{=}\mathrm{E}_{x}[\sum_{n=0}^{\tau_{B}}\mathrm{I}\left(  X_{n}=y\right)  ]$ is
the Green function of the random walk killed as it exits $B$. For $m\geq0$, we
write $G^{m}(x,y)\overset{\mathrm{def}}{=}\mathrm{E}_{x}[\sum_{n=0}
^{m}\mathrm{I}\left(  X_{n}=y\right)  ]$ for the $m$-step Green function.

Let $\mathcal{Q}$ be the covariance matrix of $p$. We write $\Vert
x\Vert_{\mathcal{Q}}=\sqrt{\left(  x,\mathcal{Q}^{-1}x\right)  },$ where
$\left(  \cdot,\cdot\right)  $ is the inner product in $\mathbb{R}^{2}.$
Observe that there exist $c^{\prime}>0$ and $c^{\prime\prime}<\infty$ such
that $c^{\prime}\left|  x\right|  \leq\Vert x\Vert_{\mathcal{Q}}\leq
c^{\prime\prime}\left|  x\right|  $.

\begin{proposition}
1. There exists a constant $K>0$ depending on $p(\,\cdot\,)$ such that
\begin{equation}
\lim_{\left|  x\right|  \rightarrow\infty}[a(x)-(\pi\sqrt{\det Q})^{-1}
\log\Vert x\Vert_{\mathcal{Q}}-K]=0\,.\label{eq_asympt_pk}%
\end{equation}

2.Let $B$ be the box of radius $R$ centered at the origin. Then, as
$R\rightarrow\infty$,
\begin{equation}
G_{B}(0,0)=(\pi\sqrt{\det Q})^{-1}\,\log R+O(1)\,.\label{eq_G_in_B}%
\end{equation}

3. Let $x\in\mathbb{Z}^{2}$. Then, as $\left|  x\right|  \rightarrow\infty$,
\begin{equation}
G_{\mathbb{Z}^{2}\setminus\{x\}}(0,0)=2(\pi\sqrt{\det Q})^{-1}\,\log\left|
x\right|  +O(1)\,.\label{eq_GF_one_obstacle}%
\end{equation}

4. As $n\rightarrow\infty$,
\begin{equation}
G^{n}(0,0)=(2\pi\sqrt{\det Q})^{-1}\log n+O(1)\,.\label{eq_G_upto_n}%
\end{equation}

5. Let $B$ be the box of radius $R$ centered at the origin, and let $x\in B$
be such that $\left|  x\right|  \leq\tfrac{1}{2}R$. Then there exist $K_{3}>0$
and $R_{0}>0$ such that, for all $R\geq R_{0}$,
\begin{equation}
\mathrm{P}_{x}[T_{\{0\}}\leq\tau_{B}]\geq\frac{K_{3}}{\log R}
\,.\label{eq_hitting}%
\end{equation}
\end{proposition}

\begin{proof}
(\ref{eq_asympt_pk}) is proved in \cite{FU}.

(\ref{eq_G_in_B}) follows from (\ref{eq_asympt_pk}) by a standard argument,
see \cite{Lawler}. The proof there is for the nearest neighbor random walk
only, but it can be easily adapted to cover the more general case considered here.

(\ref{eq_GF_one_obstacle}) follows from (\ref{eq_asympt_pk}) and P11.6 in
\cite{Spitzer}.

(\ref{eq_G_upto_n}) follows from a standard local limit theorem:
\begin{equation}
p_{n}(0)=\frac{1}{2\pi\sqrt{\det Q}n}+O(n^{-1-\varepsilon}
)\label{Est_Berry_Esseen}%
\end{equation}
for some positive $\varepsilon.$ Under the assumptions of the existence of a
third moment, this is a standard Berry-Esseen type estimate (with
$\varepsilon=1/2).$ We don't know of an exact reference under the assumption
of a $(2+\delta)$-moment only. The paper \cite{Ibra} treats the case of a
one-dimensional random walk. The method there can easily be adapted to prove
(\ref{Est_Berry_Esseen}) on the two-dimensional lattice.

Finally, (\ref{eq_hitting}) is proved in \cite{Lawler}, Proposition 1.6.7, for
the simple random walk. Again, the proof can easily been adapted to cover the
more general case.
\end{proof}

\subsection{Approximations for $p_{n}(x)$}

We will need some essentially well-known facts about $p_{n}(x)$ for large $n$
and $x,$ in case there exists an exponential moment of $p.$ For the
convenience of the reader, we sketch the argument, which is completely
standard. The results in this subsection hold for general dimensions.

\begin{proposition}
\label{PropSaddle}Assume
\[
\sum_{x}p(x)\,\mathrm{e}^{a\left|  x\right|  }<\infty
\]
for some $a>0.$ Then there exists $\eta>0$ such that for $\left|  x/n\right|
<\eta,$
\[
p_{n}(x)=\left(  \frac{1}{\left(  2\pi n\right)  ^{d/2}\sqrt{\det
\mathcal{Q}\left(  x/n\right)  }}+O\left(  \frac{1}{n^{\left(  d+1\right)
/2}}\right)  \right)  \exp\left[  -nI(x/n)\right]  ,
\]
where $\mathcal{Q}\left(  \xi\right)  ,$ $\left|  \xi\right|  <\eta,$ are
$d\times d$-matrices, depending analytically on $\xi,$ and satisfying
$\mathcal{Q}\left(  0\right)  =\mathcal{Q}.$ $I(\xi),$ $\left|  \xi\right|
<\eta,$ also depends analytically on $\xi$ and satisfies $I(0)=0,$ $\nabla
I(0)=0,$ $\nabla^{2}I(0)=\mathcal{Q}^{-1}.$
\end{proposition}

\begin{proof}
We use the standard approximation of $p_{n}(x)$ by tilting the measure and
applying a local central limit theorem with error estimate. For $\lambda
\in\mathbb{R}^{d}$ in a neighborhood of $0,$ we consider the tilted measure
\[
p^{(\lambda)}(x)\overset{\mathrm{def}}{=}\frac{p(x)\exp\left(  \lambda
,x\right)  }{z\left(  \lambda\right)  },
\]
where $z\left(  \lambda\right)  \overset{\mathrm{def}}{=}\sum_{x}%
p(x)\exp\left(  \lambda,x\right)  .$ Clearly $\nabla\log z\left(  0\right)
=0,$ and $\nabla^{2}\log z(0)=\mathcal{Q}.$ Therefore, the mapping
$\lambda\rightarrow\nabla\log z\left(  \lambda\right)  $ is an analytic
diffeomorphism of a neighborhood of $0$ to a neighborhood of $0,$ leaving $0$
fixed. Therefore, for any $\xi$ in a neighborhood of $0$ in $\mathbb{R}^{d}$,
there exists a unique $\lambda(\xi)$ with $\nabla\log z\left(  \lambda
(\xi)\right)  =\xi.$ Using this, we see that for $\left|  x\right|  \leq\eta
n,$ $\eta>0$ small enough, we can write
\[
p_{n}(x)=\exp\left[  -nI(x/n)\right]  \,p_{n}^{(\lambda(x/n))}(x),
\]
where $I(\xi)\overset{\mathrm{def}}{=}\left(  \lambda(\xi),\xi\right)  -\log
z\left(  \lambda\left(  \xi\right)  \right)  .$ Evidently, $I(0)=0,$ $\nabla
I(0)=0,$ and a simple computation yields $\nabla^{2}I(0)=\mathcal{Q}^{-1}.$
Furthermore, $p^{(\lambda(x/n))}$ has now mean exactly $x/n$ and covariance
matrix $\mathcal{Q}\left(  x/n\right)  ,$ where $\mathcal{Q}\left(
\xi\right)  $ depends analytically in $\xi$ and satisfies $\mathcal{Q}\left(
0\right)  =\mathcal{Q}.$ Applying a local central limit theorem with standard
Berry-Esseen type error estimate we get
\[
\left|  p_{n}^{(\lambda(x/n))}(x)-\frac{1}{\left(  2\pi n\right)
^{d/2}\,\sqrt{\det\mathcal{Q}\left(  x/n\right)  }}\right|  \leq\frac
{C}{n^{\left(  d+1\right)  /2}}.
\]
\end{proof}

\begin{corollary}
Assume (\ref{CondH2}). There exist $\kappa_{0}$ and $K>0$ such that for
$\kappa\geq\kappa_{0}$ and all $x\in\mathbb{Z}^{d}$ with $\left|  x\right|  $
large enough,
\begin{equation}
\mathrm{P}_{0}[X_{[\kappa\left|  x\right|  ]}=x]\geq e^{-K\,\kappa
^{-1}\,\left|  x\right|  }\,.\label{eq_hit_time}%
\end{equation}
\end{corollary}

\begin{proof}
Approximate $I$ in Proposition \ref{PropSaddle} by an appropriate quadratic function.
\end{proof}

\subsection{Crossing probabilities for thick shells}

We start with some one-dimensional considerations. Let $\left(  X_{i}\right)
_{i\geq0}$ be a $\mathbb{Z}$-valued random walk, where the distribution of the
jumps $X_{i}-X_{i-1}$ is distributed according to $\left(  q\left(  x\right)
\right)  _{x\in\mathbb{Z}},$ where we assume that $\sum_{x}x\,q(x)=0 $ and
$\sum_{x}\exp\left[  \alpha\left|  x\right|  \right]  \,q\left(  x\right)
<\infty$ for some $\alpha>0.$ We define the ladder-epochs and ladder heights
\[
\tau_{0}\overset{\mathrm{def}}{=}0,\,\xi_{0}\overset{\mathrm{def}}{=}0,
\]%
\[
\tau_{k+1}\overset{\mathrm{def}}{=}\inf\left\{  n>\tau_{k}:X_{n}\geq\xi
_{k}+1\right\}  ,\,\xi_{k+1}\overset{\mathrm{def}}{=}X_{\tau_{k+1}}.
\]
By the Markov property, the sequence $\left(  \xi_{k}-\xi_{k-1}\right)
_{k\geq1}$ is i.i.d.

\begin{lemma}
\label{LeOnedim}a)
\[
\mathrm{E}_{0}\left(  \exp\left[  \alpha^{\prime}\xi_{1}\right]  \right)
<\infty
\]
for some $\alpha^{\prime}>0.$

b) Let $K,n\in\mathbb{N},$ and define the intervals $I_{j}\overset
{\mathrm{def}}{=}((j-1)K,jK]\subset\mathbb{N}.$ Let also $\zeta\overset
{\mathrm{def}}{=}\#\left\{  j\leq n:\xi_{k}\in I_{j} ,\,\forall k\right\}  .$
Then for any $0<s<1$
\[
\limsup_{K\rightarrow\infty}\limsup_{n\rightarrow\infty}\frac{1}{Kn}
\log\mathrm{P}_{0}\left(  \zeta\leq sn\right)  <0.
\]
\end{lemma}

\begin{proof}
a) is well known. For the convenience of the reader we give a crude proof,
sufficient for our purpose.
\begin{align*}
\mathrm{P}_{0}\left(  \xi_{1}\geq k\right)   &  \leq\mathrm{P}_{0}\left(
\xi_{1}\geq k,\tau_{1}\leq\exp\left[  \lambda k\right]  \right)
+\mathrm{P}_{0}\left(  \tau_{1}>\exp\left[  \lambda k\right]  \right) \\
&  \leq\exp\left[  \lambda k\right]  \mathrm{P}_{0}\left(  X_{1}\geq k\right)
+\frac{C}{\sqrt{\exp\left[  \lambda k\right]  }}\leq\exp\left[  -\alpha
^{\prime}k\right]
\end{align*}
for some $\alpha^{\prime}>0,$ by choosing $\lambda>0$ appropriately, for large
enough $k.$

b) Let $\sigma\overset{\mathrm{def}}{=}\min\left\{  j:\xi_{j}>Kn\right\}  .$
Then, by standard large deviation estimates,
\[
\mathrm{P}_{0}\left(  \sigma>\lambda Kn\right)  =\mathrm{P}_{0}\left(
\xi_{\lambda Kn}<Kn\right)  \leq\exp\left[  -CKn\right]
\]
for $Kn$ large enough, when $\lambda$ is chosen appropriately (e.g.
$\lambda=1/2\mathrm{E}_{0}\xi_{1}).$ We consider the independent differences
$\Delta_{j}\overset{\mathrm{def}}{=}\xi_{j}-\xi_{j-1}.$ We have for $0<s<1:$
\begin{align*}
\mathrm{P}_{0}\left(  \zeta\leq sn,\sigma\leq\lambda Kn\right)   &
\leq\mathrm{P}_{0}\left(  \sum\nolimits_{j=1}^{\lambda Kn}\frac{\Delta_{j}}%
{K}\mathrm{I}\left(  \Delta_{j}>K\right)  \geq(1-s)n\right) \\
&  \leq\exp\left[  -aK(1-s)n\right]  \left\{  \mathrm{E}_{0}\left(
\exp\left[  a\Delta_{j}\mathrm{I}\left(  \Delta_{j}>K\right)  \right]
\right)  \right\}  ^{\lambda Kn},
\end{align*}
for any $a>0.$ According to a), we can choose $a$ such that $\mathrm{E}%
_{0}\left(  \exp\left[  a\Delta_{j}\right]  \right)  <\infty,$ and then, for
any $\delta>0,$ we may choose $K$ large enough, such that $\mathrm{E}%
_{0}\left(  \exp\left[  a\Delta_{j}\mathrm{I}\left(  \Delta_{j}>K\right)
\right]  \right)  \leq1+\delta,$ i.e.
\[
\mathrm{P}_{0}\left(  \zeta\leq sn,\sigma\leq\lambda Kn\right)  \leq
\exp\left[  -aK(1-s)n+\delta\lambda Kn\right]  \leq\exp\left[  \frac
{-aK(1-s)n}{2}\right]  ,
\]
if $\delta\leq aK(1-s)/2.$ This proves the claim.
\end{proof}

We apply this now to our $d$-dimensional random walk, where we again assume
the existence of an exponential moment (\ref{CondH2}). Let $\mathbb{Z}^{d}(K)
$ be the division of $\mathbb{Z}^{d}$ into square cells of side length $K,$
where we take $\left\{  1,\ldots,K\right\}  ^{d}$ as one of the block. We
further consider a big square $S_{n,K}\overset{\mathrm{def}}{=}\left\{
-nK+1,-nK+2,\ldots,nK\right\}  ^{d},$ which of course is divided in $\left(
2n\right)  ^{d}$ cells of side-length $K.$

\begin{proposition}
\label{PropCross}Let $\eta_{n,K}$ be the number of cells in $\mathbb{Z}%
^{d}(K)$ which are visited by the random walk up to time $\tau_{S_{n,K}}.$
Then for any $s\in\left(  0,1\right)  $
\[
\limsup_{K\rightarrow\infty}\limsup_{n\rightarrow\infty}\frac{1}{Kn}%
\log\mathrm{P}_{0}\left(  \eta_{n,K}\leq sn\right)  <0.
\]
\end{proposition}

\begin{proof}
The proposition is an easy consequence of the one-dimensional result. Indeed,
write the random walk in (dependent) components $X_{n}=\left(  X_{n,1}%
,X_{n,2},\ldots,X_{n,d}\right)  ,$ where $\left(  X_{n,i}\right)  $ are
one-dimension random walks, possessing an exponential moment. The first time
$\tau_{S_{n,K}}$ when $\left(  X_{n}\right)  $ leaves $S_{n,K}$ is also the
first time where one of the $\left(  X_{n,i}\right)  $ $1\leq i\leq d$ leaves
the interval $\left\{  -nK+1,-nK+2,\ldots,nK\right\}  .$ Assume for instance
that at $\tau_{S_{n,K}}$, $\left(  X_{n,1}\right)  $ for the first time leaves
the above interval to the right. (There are of course $2d-1$ other cases).
This is then the first time it surpasses $nK.$ Furthermore, the number of
$K$-cells visited by the $d$-dimensional walk is at least the number of
intervals among $(1,K],\,(K,2K],\,\ldots\,,((n-1)K,nK],$ visited by $\left(
X_{n,1}\right)  .$ For the other $2d-1$ cases, similar statements hold, of
course. From this observation, Proposition \ref{PropCross} follows immediately
from Lemma \ref{LeOnedim}.
\end{proof}

\section{On the range of a random walk\label{RangeRW}}

We present here two results about the number of points visitied by a
two-dimensional random walk.

\subsection{Tied-down expectations of $\left|  X_{[0,n]}\right|  $}

We write \textrm{$P$}$_{x,y}^{(n)}$ for the law of a random walk, starting in
$x,$ conditioned on $X_{n}=y.$ Of course, we tacitly always assume that the
probability of the latter is positive whenever we use this notation, which is
certainly true for large enough $n$ (depending on $x,y).$ We will need some
information on the first return probabilities:
\[
q_{l}\overset{\mathrm{def}}{=}\mathrm{P}_{0}\left(  X_{1}\neq0,\ldots
,X_{l-1}\neq0,X_{l}=0\right)  .
\]
By recurrence of the random walk, we have $\sum_{l\geq1}q_{l}=1,$ and the
following estimate is well-known.

\begin{lemma}
\label{LeTaboo}
\begin{equation}
q_{l}=\frac{\gamma}{l\left(  \log l\right)  ^{2}}+o\left(  \frac{1}{l\left(
\log l\right)  ^{2}}\right)  ,\label{Taboo}%
\end{equation}
where $\gamma>0,$ as $l\rightarrow\infty.$
\end{lemma}

We need some information on $\mathrm{E}_{0,x}^{(n)}\left(  \left|
X_{[0,n]}\right|  \right)  .$

\begin{proposition}
\label{PropTieddown}There exist $A_{0}>1$ and $C>0$ such that for all $A\geq
A_{0}$, there exists $r_{0}(A)\in\mathbb{N},$ such that for $\left|  x\right|
\geq r_{0}(A)$ and $n$ defined by
\[
n\overset{\mathrm{def}}{=}\left[  A\left|  x\right|  \right]  ,
\]
one has
\[
\mathrm{E}_{0,x}^{(n)}\left(  \left|  X_{[0,n]}\right|  \right)  \leq
C\frac{n}{\log A}.
\]
\end{proposition}

\begin{proof}
Remark first that under the conditions of the lemma, $\mathrm{P}_{0}\left(
X_{n}=x\right)  >0$ for the large enough $r_{0}(A).$ This easily follows from
irreducibility and aperiodicity. Therefore, $\mathrm{E}_{0,x}^{(n)}\left(
\left|  X_{[0,n]}\right|  \right)  $ is well-defined. We first derive a simple
exact expression for this expectation:
\begin{equation}
\mathrm{E}_{0,x}^{(n)}\left(  \left|  X_{[0,n]}\right|  \right)
=n+1-\sum_{l=1}^{n}\left(  n-l+1\right)  q_{l}\frac{p_{n-l}(x)}{p_{n}
(x)}.\label{eqTiedsaus}%
\end{equation}
This readily follows from a standard ``last exit - first entrance'' decomposition.%

\[
\mathrm{E}_{0,x}^{(n)}\left(  \left|  X_{[0,n]}\right|  \right)
=\frac{\mathrm{E}_{0}\left(  \left|  X_{[0,n]}\right|  ;X_{n}=x\right)
}{p_{n}(x)}.
\]%
\begin{align}
\mathrm{E}_{0}\left(  \left|  X_{[0,n]}\right|  ;X_{n}=x\right)   &
=\sum_{y\in\mathbb{Z}^{2}}\mathrm{P}_{0}\left(  X_{k}=y\,\mathrm{for\,some\,}
k\in\lbrack0,n],X_{n}=x\right) \nonumber\\
&  =\sum_{y\in\mathbb{Z}^{2}}\sum_{k=0}^{n}\mathrm{P}_{0}\left(  X_{k}
=y,X_{k+1}\neq y,\ldots,X_{n-1}\neq y,X_{n}=x\right) \label{Split1}\\
&  =\sum_{y\in\mathbb{Z}^{2}}\sum_{k=0}^{n}\mathrm{P}_{0}\left(  X_{k}
=y\right)  \mathrm{P}_{y}\left(  X_{1}\neq y,\ldots,X_{n-k-1}\neq
y,X_{n-k}=x\right)  .\nonumber
\end{align}%
\begin{align*}
\mathrm{P}_{y}\bigl(  X_{1}\neq y,\ldots,&X_{n-k-1}\neq y,X_{n-k}=x\bigr)  \\ 
&=p_{n-k}\left(  x-y\right)  -\sum_{l=1}^{n-k}\mathrm{P}_{y}\left(  X_{1}\neq
y,\ldots,X_{l-1}\neq y,X_{l}=y,X_{n-k}=x\right) \\
&  =p_{n-k}\left(  x-y\right)  -\sum_{l=1}^{n-k}q_{l}\,p_{n-k-l}(x-y).
\end{align*}
Implementing this into (\ref{Split1}) and summing over $y$ yields
\begin{align*}
E_{0}\left(  \left|  X_{[0,n]}\right|  ;X_{n}=x\right)   &  =\left(
n+1\right)  p_{n}(x)-\sum_{k=0}^{n-1}\sum_{l=1}^{n-k}q_{l}\,p_{n-l}(x)\\
&  =\left(  n+1\right)  p_{n}(x)-\sum_{l=1}^{n}\left(  n-l+1\right)
q_{l}\,p_{n-l}(x).
\end{align*}
From this, (\ref{eqTiedsaus}) follows.

We next use this together with the information on $q_{l}$ in Lemma
\ref{LeTaboo} and Proposition \ref{PropSaddle} to get the desired estimate.%

\begin{align}
E_{0,x}^{(n)}\left(  \left|  X_{[0,n]}\right|  \right)   &  \leq
n+1-\sum_{l=1}^{A}\left(  n-l+1\right)  q_{l}\frac{p_{n-l}(x)}{p_{n}
(x)}\nonumber\\
&  =\left(  n+1\right)  \left[  1-\sum_{l=1}^{A}q_{l}\,\frac{\left(
n-l+1\right)  p_{n-l}(x)}{\left(  n+1\right)  p_{n}(x)}\right] \nonumber\\
&  =\left(  n+1\right)  \left[  1-\sum_{l=1}^{A}q_{l}+\sum_{l=1}^{A}
q_{l}\left(  1-\frac{\left(  n-l+1\right)  p_{n-l}(x)}{\left(  n+1\right)
p_{n}(x)}\right)  \right] \label{Est1}\\
&  =\left(  n+1\right)  \sum_{l=A+1}^{\infty}q_{l}+\left(  n+1\right)
\sum_{l=1}^{A}q_{l}\left(  1-\frac{\left(  n-l+1\right)  p_{n-l}(x)}{\left(
n+1\right)  p_{n}(x)}\right)  ,\nonumber
\end{align}
the last equation by recurrence of the two-dimensional random walk. From Lemma
\ref{LeTaboo}, we get
\begin{equation}
\sum_{l=A+1}^{\infty}q_{l}\leq\frac{C}{\log A},\label{EstTailQ}%
\end{equation}
and it therefore suffices to estimate the second summand on the right-hand
side of (\ref{Est1}). If $n$ is large enough (depending on $A),$ then $\left|
x\right|  \leq2\left(  n-l\right)  /A$ whenever $l\leq A.$ We use Proposition
\ref{PropSaddle} and obtain for $\left|  x\right|  \rightarrow\infty,$ and
therefore $n\rightarrow\infty$
\[
\frac{\left(  n-l-1\right)  p_{n-l}(x)\mathrm{\exp}\left[  -\left(
n-l\right)  I\left(  \frac{x}{n-l}\right)  \right]  }{\left(  n+1\right)
p_{n}(x)\mathrm{\exp}\left[  -nI\left(  \frac{x}{n}\right)  \right]  }
=1+O_{A}\left(  \frac{1}{\sqrt{n}}\right)  ,
\]
where $O_{A}\left(  \frac{1}{\sqrt{n}}\right)  $ means that there is a
constant $c_{A}$ depending on $A$ such that $\left|  O_{A}\left(  \frac
{1}{\sqrt{n}}\right)  \right|  \leq\frac{c_{A}}{\sqrt{n}}.$ Furthermore, by
Taylor expansion, we get
\begin{align*}
&  \exp\left[  -\left(  n-l\right)  I\left(  \frac{x}{n-l}\right)  +nI\left(
\frac{x}{n}\right)  \right] \\
&  =\exp\left[  l\left\{  I\left(  \frac{x}{n}\right)  -\left(  \frac{x}
{n},\nabla I\left(  \frac{x}{n}\right)  \right)  \right\}  +O_{A}\left(
\frac{1}{n}\right)  \right]  .
\end{align*}
Remark that
\[
\left|  I\left(  \frac{x}{n}\right)  -\left(  \frac{x}{n},\nabla I\left(
\frac{x}{n}\right)  \right)  \right|  \leq C\left(  \frac{\left|  x\right|
}{n}\right)  ^{2}.
\]
Combining these observations, we get
\begin{align*}
\sum_{l=1}^{A}q_{l}\left|  1-\frac{\left(  n-l+1\right)  p_{n-l}(x)}{\left(
n+1\right)  p_{n}(x)}\right|   &  \leq\sum_{l=1}^{A}q_{l}\left(  \exp\left[
C\frac{l}{A^{2}}\right]  -1\right)  +O_{A}\left(  \frac{1}{\sqrt{n}}\right) \\
&  \leq\frac{C}{A^{2}}\sum_{l=1}^{A}\frac{1}{\left(  \log l\right)  ^{2}}
+O_{A}\left(  \frac{1}{\sqrt{n}}\right) \\
&  \leq\frac{C}{A}+O_{A}\left(  \frac{1}{\sqrt{n}}\right)  \leq\frac{2C}{A},
\end{align*}
for large enough $n$ ($n\geq n_{0}\left(  A\right)  ).$ This is much better
than required, and therefore proves the proposition.
\end{proof}

\subsection{Moderate deviations for $\left|  X_{[0,n]}\right|  $}

We use a variant of the approach in \cite{vdBBdH} to prove the following result.

\begin{proposition}
\label{PropModerate}Assume (\ref{CondH1}). For any $R>0$ there exists
$\kappa>0$ such that
\begin{equation}
\mathrm{P}_{0}\left(  \left|  X_{[0,n]}\right|  \leq\kappa\,\frac{n}{\log
n}\right)  \leq n^{-R}\,,
\end{equation}
for all $n$ large enough.
\end{proposition}

In contrast to our standard convention about constants denoted by $C$,
$c_{1},c_{2},\ldots$ are positive constants which are always the same after
they had been introduced. If these constants depend on other parameters, it
will be clearly indicated. All inequalities are supposed to hold only for
large enough $n$ without further notice.

Let $L_{n}\overset{\mathrm{def}}{=}\left[  \sqrt{n/\log n}\right]  $ and
$T_{n}\overset{\mathrm{def}}{=}\{0,1,\ldots,L_{n}-1\}^{d}.$ We periodize the
random walk by setting
\[
\hat{X}_{n}\overset{\mathrm{def}}{=}X_{n}\,\operatorname{mod}\,L_{n},
\]
coordinatewise, getting therefore a random walk on the discrete torus $T_{n}.
$ The transition probabilities are given by $\hat{p}(x)\overset{\mathrm{def}
}{=}\sum_{y=x\,\operatorname{mod}\,L_{n}}p(y).$ The number of points $\bigl
|\hat{X}_{[0,n]}\bigr|$ visited by the periodized random walk is clearly at
most $\bigl|X_{[0,n]}\bigr|$. Therefore
\[
\mathrm{P}_{0}\left(  \bigl|X_{[0,n]}\bigr|\leq\kappa\frac{n}{\log n}\right)
\leq\mathrm{P}_{0}\left(  \bigl|\hat{X}_{[0,n]}\bigr|\leq\kappa\frac{n}{\log
n}\right)  .
\]
For the rest of this section, we always work with this periodized walk, but
leave the hat $\symbol{94}$ out in the notations for the sake of notational
convenience. For convenience, we also assume that $p$ is aperiodic. The
general case requires only some trivial adjustments.

We choose $m=m_{n}=\left[  \delta\frac{n}{\log n}\right]  ,$ where $\delta>0$
is a (small) number, to be specified later on. We also set $K=K_{n}=\left[
n/m_{n}\right]  \approx\frac{\log n}{\delta}.$ We denote by $\mathbf{X}$ the
sequence of points observed at multiples of $m:$
\[
\mathbf{X}\overset{\mathrm{def}}{=}\left(  X_{0},X_{m},X_{2m},\ldots
,X_{Km}\right)  .
\]
The set of points (on the torus) visited during the $i-$th time interval is
denoted by $V_{i}^{0}:$
\[
V_{i}^{0}\overset{\mathrm{def}}{=}\left\{  X_{(i-1)m+1},X_{(i-1)m+2}
,\ldots,X_{im}\right\}  .
\]
We introduce a truncation by defining
\[
V_{i}\overset{\mathrm{def}}{=}\left\{
\begin{array}
[c]{cc}%
V_{i}^{0} & \mathrm{if\,}d(X_{(i-1)m},X_{im})\leq b\sqrt{m}\\
\emptyset & \mathrm{if\,}d(X_{(i-1)m},X_{im})>b\sqrt{m}%
\end{array}
\right.  .
\]
$d$ is the lattice distance on the discrete torus. We also write
$\mathbf{V}\overset{\mathrm{def}}{=}(V_{1},V_{2},\ldots,V_{K}).$ Remark that
$d(X_{(i-1)m},X_{im})$ are i.i.d. random variables and
\[
\mathrm{P}_{0}\left(  d(X_{(i-1)m},X_{im})>b\sqrt{m}\right)  \leq\exp\left[
-c_{1}b^{2}\right]  .
\]
Let
\[
\Gamma_{n,b,\delta}\overset{\mathrm{def}}{=}\#\left\{  i:d(X_{(i-1)m}
,X_{im})>b\sqrt{m}\right\}  .
\]
Then $\Gamma_{n,b,\delta}$ is binomially distributed, and we obtain
\begin{equation}
\mathrm{P}_{0}\left(  \Gamma_{n,b,\delta}\geq2\exp\left[  -c_{1}b^{2}\right]
\frac{\log n}{\delta}\right)  \leq\exp\left[  -c_{2}\frac{\log n}{\delta}
\exp\left[  -c_{1}b^{2}\right]  \right] \label{EstGamma}%
\end{equation}

We denote by$\mathcal{\ P}(T_{n})$ the set of subsets of $T_{n}$ and by
$\Psi:\mathcal{P}(T_{n})^{K}\rightarrow\mathbb{N}$ the mapping
\[
\mathbf{V}\rightarrow\left|  \bigcup\nolimits_{i=1}^{K}V_{i}\right|  .
\]
Clearly, $\Psi$ is Lipschitz in the sense
\[
\left|  \Psi\left(  \mathbf{V}\right)  -\Psi\left(  \mathbf{U}\right)
\right|  \leq\sum_{i=1}^{K}\left|  V_{i}\bigtriangleup U_{i}\right|  .
\]
Using this notation, we get
\begin{align*}
\mathrm{P}_{0}\left(  \left| X_{[0,n]}\right| \leq\kappa\frac{n}{\log
n}\right)   &  \leq\mathrm{P}_{0}\left(  \Psi\left(  \mathbf{V}\right)
\leq\kappa\frac{n}{\log n}\right) \\
&  =\mathrm{E}_{0}\left(  \mathrm{P}_{\mathbf{X}}\left(  \Psi\left(
\mathbf{V}\right)  \leq\kappa\frac{n}{\log n}\right)  \right)  ,
\end{align*}
where $\mathrm{P}_{\mathbf{X}}$ denotes the conditional law given the vector
$\mathbf{X}.$ Under $\mathrm{P}_{\mathbf{X}}$, the sets $V_{i}$ are
independent random subsets of the torus $T_{n}$. We thus can apply a general
result of Talagrand. Let $\mu=\mu_{\mathbf{X}}$ be a median of the
(conditional) distribution of $\Psi,$ i.e. a number with $\mathrm{P}
_{\mathbf{X}}\left(  \Psi\left(  \mathbf{V}\right)  \leq\mu_{\mathbf{X}
}\right)  \geq1/2$ and $\mathrm{P}_{\mathbf{X}}\left(  \Psi\left(
\mathbf{V}\right)  \geq\mu_{\mathbf{X}}\right)  \geq1/2.$ Let $f:
\mathcal{\text{P}}(T_{n})^{K}\rightarrow\mathbb{N}$ be defined by
\[
f\left(  \mathbf{V}\right)  \overset{\mathrm{def}}{=}\inf\left\{  \sum
_{i=1}^{K}\left|  V_{i}\bigtriangleup U_{i}\right|  :\Psi\left(
\mathbf{U}\right)  \leq\mu_{\mathbf{X}}\right\}  .
\]
Then by Theorem 2.4.1 of \cite{Tala}, we have for any $a>0$ and $\lambda>0$
\[
\mathrm{P}_{\mathbf{X}}\left(  f\left(  \mathbf{V}\right)  \geq a\right)
\leq\Xi_{\mathbf{X}}(a,\lambda),
\]
where
\[
\Xi_{\mathbf{X}}(a,\lambda)\overset{\mathrm{def}}{=}2\mathrm{e}^{-\lambda
a}\prod_{i=1}^{K}\mathrm{E}_{\mathbf{X}}\left(  \cosh\left(  \lambda\left|
V_{i}\bigtriangleup U_{i}\right|  \right)  \right)  ,
\]
and where the $\mathbf{U}$ is an independent copy of $\mathbf{V}$ (under the
conditional law). Similarly, putting
\[
\widehat{f}\left(  \mathbf{V}\right)  \overset{\mathrm{def}}{=}\inf\left\{
\sum_{i=1}^{K}\left|  V_{i}\bigtriangleup U_{i}\right|  :\Psi\left(
\mathbf{U}\right)  \geq\mu_{\mathbf{X}}\right\}  ,
\]
we get
\[
\mathrm{P}_{\mathbf{X}}\left(  \widehat{f}\left(  \mathbf{V}\right)  \geq
a\right)  \leq\Xi_{\mathbf{X}}(a,\lambda).
\]
Combining these two estimates, we get
\[
\mathrm{P}_{\mathbf{X}}\left(  \left|  \Psi\left(  \mathbf{V}\right)
-\mu_{\mathbf{X}}\right|  \geq a\right)  \leq2\Xi_{\mathbf{X}}(a,\lambda).
\]
Now,
\begin{align*}
\mathrm{P}_{\mathbf{X}}\left(  \Psi\left(  \mathbf{V}\right)  \leq a\right)
&  \leq\mathrm{P}_{\mathbf{X}}\left(  \left|  \Psi\left(  \mathbf{V}\right)
-\mu_{\mathbf{X}}\right|  \geq a\right) \\
&  +\mathrm{I}\left(  \left|  \mu_{\mathbf{X}}-\mathrm{E}_{\mathbf{X}}%
\Psi\left(  \mathbf{V}\right)  \right|  \geq2a\right)  +\mathrm{I}\left(
\mathrm{E}_{\mathbf{X} }\Psi\left(  \mathbf{V}\right)  \leq4a\right)  .
\end{align*}
Remark that
\[
\left|  \mu_{\mathbf{X}}-\mathrm{E}_{\mathbf{X}}\Psi\left(  \mathbf{V}\right)
\right|  \leq a+\left|  T_{n}\right|  \mathrm{P}_{\mathbf{X}}\left(  \left|
\Psi\left(  \mathbf{V}\right)  -\mu_{\mathbf{X}}\right|  \geq a\right)  ,
\]
and therefore
\begin{align}
\mathrm{P}_{\mathbf{X}}\left(  \Psi\left(  \mathbf{V}\right)  \leq a\right)
&  \leq\mathrm{P}_{\mathbf{X}}\left(  \left|  \Psi\left(  \mathbf{V}\right)
-\mu_{\mathbf{X}}\right|  \geq a\right) \nonumber\\
&  +\mathrm{I}\left(  \mathrm{P}_{\mathbf{X}}\left(  \left|  \Psi\left(
\mathbf{V} \right)  -\mu_{\mathbf{X}}\right|  \geq a\right)  \geq\frac
{a}{\left|  T_{n}\right|  }\right)  +\mathrm{I}\left(  \mathrm{E}_{\mathbf{X}%
}\Psi\left(  \mathbf{V}\right)  \leq4a\right) \label{EstPsi}\\
&  \leq2\Xi_{\mathbf{X}}(a,\lambda)+\mathrm{I}\left(  \Xi_{\mathbf{X}%
}(a,\lambda)\geq\frac{a}{2\left|  T_{n}\right|  }\right)  +\mathrm{I}\left(
\mathrm{E} _{\mathbf{X}}\Psi\left(  \mathbf{V}\right)  \leq4a\right)
.\nonumber
\end{align}
We apply this inequality to $a=a_{n}\overset{\mathrm{def}}{=}\kappa\frac
{n}{\log n},$ and with $\lambda=\lambda_{n}\overset{\mathrm{def}}{=}%
A\frac{\left(  \log n\right)  ^{2}}{n},$ where $A$ will be specified below.
Then we have
\[
\Xi_{\mathbf{X}}\left(  \kappa\frac{n}{\log n},A\frac{\left(  \log n\right)
^{2}}{n}\right)  =2\exp\left[  -A\kappa\log n\right]  \prod_{i=1}%
^{K}\mathrm{E}_{\mathbf{X}}\left(  \cosh\left(  A\frac{\left(  \log n\right)
^{2}}{n}\left|  V_{i}\bigtriangleup U_{i}\right|  \right)  \right)  .
\]%
\[
\mathrm{E}_{\mathbf{X}}\left(  \cosh\left(  A\frac{\left(  \log n\right)
^{2}}{n}\left|  V_{i}\bigtriangleup U_{i}\right|  \right)  \right)
\leq\mathrm{E}_{\mathbf{X}}\left(  \cosh\left(  2A\delta\frac{\log m}
{m}\left|  V_{i}\bigtriangleup U_{i}\right|  \right)  \right)  .
\]
We assume now
\begin{equation}
2A\delta<1,\label{CondA1}%
\end{equation}
and use $\cosh(xy)\leq1+x^{2}\mathrm{e}^{y}$ for $0\leq x\leq1,$ $0\leq y. $
Furthermore, we have the following

\begin{lemma}
\label{LeExponMoment}
\[
\mathrm{E}_{\mathbf{X}}\exp\left[  \frac{\log m}{m}V_{i}\right]  \leq C(b).
\]
\end{lemma}

\begin{proof}
We can take $i=1.$ If $d(0,X_{m})>b\sqrt{m},$ then $V_{1}=\emptyset$ and there
is nothing to prove.

We write $\mathrm{P}_{0,x}^{(m)}$ for the law of the random walk $(X_{0}
,X_{1},\ldots,X_{m})$ conditioned on $X_{0}=0,$ $X_{m}=x.$ (For simplicity, we
neglect trivial parity problems.) Let $Z_{T}(m/2)$ be the number of points
visited by $X_{1},\ldots,X_{m/2}$ on the torus (assuming $m $ for simplicity
to be even). Then it suffices to prove for $d(0,x)\leq b\sqrt{m}$
\begin{equation}
\mathrm{E}_{0,x}^{(m)}\left(  \exp\left[  \frac{\log m}{m}Z_{T}(m/2)\right]
\right)  \leq C(b).\label{EstSaus}%
\end{equation}
The left hand side of this equals
\[
\sum_{y}\mathrm{E}_{0}\left(  \exp\left[  \frac{\log m}{m}Z_{T}(m/2)\right]
\mathrm{I}\left( \{ X_{m/2}=y\}\right)  \}\right)  \frac{p_{m/2}(x-y)}%
{p_{m}(x)}\leq C(b)\mathrm{E}_{0}\exp\left[  \frac{\log m}{m}Z_{T}%
(m/2)\right]  ,
\]
because for $d(0,x)\leq b\sqrt{m}$ we have $p_{m}(x)\geq C(b)/m>0,$ and for
all $y,$ $p_{m/2}(x-y)\leq C^{\prime}(b)/m.$ We can replace $Z_{T}(m/2)$ by
$Z(m),$ the number of points visited by a random walk of length $m$ on
$\mathbb{Z}^{2}.$ (We replace $m/2$ by $m$ just for notational convenience).
\begin{equation}
\mathrm{E}_{0}\exp\left[  \frac{\log m}{m}Z(m)\right]  \leq\mathrm{E}_{0}%
\exp\left[  \frac{\log m}{m}Z\left(  \frac{m}{\left(  \log m\right)  ^{3}%
}\right)  \right]  ^{\left(  \log m\right)  ^{3}},\label{saus1}%
\end{equation}
by the Markov property. We write $Z^{\prime}$ for $Z\left(  \frac{m}{\left(
\log m\right)  ^{3}}\right)  $.
\begin{align}
\mathrm{E}_{0}\exp\left[  \frac{\log m}{m}Z^{\prime}\right]   &  \leq
1+\frac{\log m}{m}\mathrm{E}_{0}Z^{\prime}+\frac{\left(  \log m\right)  ^{2}%
}{2m^{2}} \mathrm{E}_{0}\left(  Z^{\prime2}\exp\left[  \frac{\log m}%
{m}Z^{\prime}\right]  \right) \label{saus2}\\
&  \leq\exp\left[  \frac{\log m}{m}\mathrm{E}_{0}Z^{\prime}+\frac{\left(  \log
m\right)  ^{2}}{2m^{2}}\mathrm{E}_{0}\left(  Z^{\prime2}\exp\left[  \frac{\log
m}{m}Z^{\prime}\right]  \right)  \right]  .\nonumber
\end{align}%
\begin{align*}
\mathrm{E}_{0}\left(  Z^{\prime2}\exp\left[  \frac{\log m}{m}Z^{\prime
}\right]  \right)   &  \leq\sqrt{\mathrm{E}_{0}Z^{\prime4}\mathrm{E}_{0}%
\exp\left[  2\frac{\log m}{m}Z^{\prime}\right]  }\\
&  \leq C\frac{m^{2}}{\left(  \log m\right)  ^{6}}\sqrt{\mathrm{E}_{0}%
\exp\left[  2\frac{\left(  \log m\right)  ^{3}}{m}Z^{\prime}\right]
\mathrm{E}_{0}\exp\left[  2\frac{\log m}{m}Z^{\prime}\right]  }\\
&  \leq C\frac{m^{2}}{\left(  \log m\right)  ^{6}}\mathrm{E}_{0}\exp\left[
2\frac{\left(  \log m\right)  ^{3}}{m}Z^{\prime}\right] \\
&  \leq C\frac{m^{2}}{\left(  \log m\right)  ^{6}}\left\{  \mathrm{E}_{0}
\exp\left[  2Z^{\prime}\right]  \right\}  ^{\frac{\left(  \log m\right)  ^{3}
}{m}}\leq C\frac{m^{2}}{\left(  \log m\right)  ^{6}}.
\end{align*}
Implementing this into (\ref{saus1}) and (\ref{saus2}), this gives
\[
\mathrm{E}_{0}\exp\left[  \frac{\log m}{m}Z(m)\right]  \leq\exp\left[
\frac{\left(  \log m\right)  ^{4}}{m}\mathrm{E}_{0}Z^{\prime}+\frac{C}{\log
m}\right]  .
\]
As
\[
\mathrm{E}_{0}Z^{\prime}\leq C\frac{\frac{m}{\left(  \log m\right)  ^{3}}}
{\log\left(  \frac{m}{\left(  \log m\right)  ^{3}}\right)  }\leq C\frac
{m}{\left(  \log m\right)  ^{4}},
\]
this proves the claim.
\end{proof}

Using this lemma, we get
\[
\mathrm{E}_{\mathbf{X}}\left(  \cosh\left(  2A\delta\frac{\log m}{m}\left|
V_{i}\bigtriangleup U_{i}\right|  \right)  \right)  \leq1+(2A\delta)^{2}C(b).
\]
Therefore, we obtain
\begin{align*}
\Xi_{\mathbf{X}}\left(  \kappa\frac{n}{\log n},A\frac{\left(  \log n\right)
^{2}}{n}\right)   &  \leq2\exp\left[  -A\kappa\log n\right] \\
&  \times\left(  1+c_{3}(b)(2A\delta)^{2}\right)  ^{(1/\delta)\log n}\\
&  \leq2\exp\left[  -\frac{A\kappa}{2}\log n\right]  ,
\end{align*}
if
\[
8c_{3}(b)A\delta<\kappa.
\]
We fix
\[
A(\kappa,\delta,b)\overset{\mathrm{def}}{=}\frac{\kappa} {16c_{3}(b)\delta}.
\]
Remark that we are then also on the safe side concerning (\ref{CondA1})
provided $\kappa\leq\kappa_{0}(b),$ $\kappa_{0}(b)$ small enough. Therefore
\[
\Xi_{\mathbf{X}}\left(  \kappa\frac{n}{\log n},A\frac{\left(  \log n\right)
^{2}}{n}\right)  \leq\exp\left[  -\frac{\kappa^{2}} {24c_{3}(b)\delta}\log
n\right]  .
\]
This is a deterministic bound. We see that the second summand on the right
hand side of (\ref{EstPsi}) is zero (with $a=\kappa\frac{n}{\log n})$ for $n$
large enough, and therefore
\begin{equation}
\mathrm{P}_{0}\left(  \Psi(\mathbf{V})\leq\kappa\frac{n}{\log n}\right)
\leq\mathrm{P}_{0}\left(  \mathrm{E}_{\mathbf{X}}\Psi(\mathbf{V})\leq
4\kappa\frac{n}{\log n}\right)  +\exp\left[  -\frac{\kappa^{2} }%
{32c_{3}(b)\delta}\log n\right]  .\label{Est2Psi}%
\end{equation}
We choose now
\begin{equation}
\delta\overset{\mathrm{def}}{=}\kappa^{3},\label{DefDelta}%
\end{equation}
and so the second summand in (\ref{Est2Psi}) is fine, again if $\kappa
\leq\kappa_{0}(b),$ $\kappa_{0}(b)$ small enough. The reader should keep in
mind that $\mathbf{V}$ depends on our truncation parameter $b,$ which we
emphasize by writing $\mathbf{V}_{b}.$ Combining what we have achieved so far,
we see that it suffices to prove that for any $R>0$ there exists $b$ (large
enough) and then $\kappa>0$ small enough (depending on $b)$ such that
\begin{equation}
\mathrm{P}_{0}\left(  \mathrm{E}_{\mathbf{X}}\Psi(\mathbf{V}_{b})\leq
\kappa\frac{n}{\log n}\right)  \leq n^{-R}.\label{Task}%
\end{equation}%

\begin{align*}
\mathrm{E}_{\mathbf{X}}\Psi(\mathbf{V}_{b})  &  =\sum_{x\in T_{n}}
\mathrm{P}_{\mathbf{X}}\left(  \bigcup\nolimits_{i=1}^{K}\left\{  x\in
V_{i,b}\right\}  \right) \\
&  =\sum_{x\in T_{n}}\left(  1-\prod\nolimits_{i=1}^{K}\left(  1-\mathrm{P}
_{\mathbf{X}}\left(  x\in V_{i,b}\right)  \right)  \right) \\
&  \geq\sum_{x\in T_{n}}\left(  1-\exp\left[  -\sum\nolimits_{i=1}
^{K}\mathrm{P}_{\mathbf{X}}\left(  x\in V_{i,b}\right)  \right]  \right)  .
\end{align*}
We now chop the torus $T_{n}$ into $M=1/\delta$ subsquares $S_{1},S_{2}
,\ldots,S_{M}$ of sidelength $\sqrt{\delta\frac{n}{\log n}}$. For notational
convenience, we will assume that $\sqrt{1/\delta}$ is an integer, which
evidently is no restriction. (Remember the setting $\delta=\kappa^{3}$ but for
the moment, this will be of no importance). We set
\[
\xi_{i}\overset{\mathrm{def}}{=}\#\left\{  j\in\left\{  1,\ldots,K\right\}
:X_{(j-1)m}\in S_{i},\,d(X_{(j-1)m},X_{jm})<b\sqrt{m}\right\}
\]
and
\[
\overline{\xi}_{i}\overset{\mathrm{def}}{=}\#\left\{  j\in\left\{
1,\ldots,K\right\}  :X_{(j-1)m}\in S_{i}\right\}
\]

\begin{lemma}
Let $X_{(j-1)m}\in S_{i},$ $d(X_{(j-1)m},X_{jm})\leq b\sqrt{m},$ and $x\in
S_{i}.$ Then
\[
\mathrm{P}_{\mathbf{X}}\left(  x\in V_{j,b}\right)  \geq\frac{c_{4}(b)}{\log
n}.
\]
\end{lemma}

\begin{proof}
We use the same notations as in the proof of Lemma \ref{LeExponMoment}:
$\mathrm{P}_{y,z}^{(m)}$ denotes the law of the random walk of length $m$ (on
the torus), conditioned to start in $y$ and to end in $z.$ If $x,y\in S_{i}$
and $d(y,z)\leq b\sqrt{m},$ then
\begin{multline*}
\mathrm{P}_{y,z}^{(m)}\left(  X_{j}=x\,\mathrm{for\,some\,}j\in\left\{
1,\ldots,m\right\}  \right)    \geq\mathrm{P}_{y,z}^{(m)}\left(
X_{j}=x\,\mathrm{for\,some\,}j\in\lbrack m/4,m/2]\right) \\
=\frac{\sum_{j=m/4}^{m/2}p_{j}(x-y)\mathrm{P}_{0}\left(  X_{1}\neq
0,\ldots,X_{m/2-j-1}\neq0,X_{m-j}=z-x\right)  }{p_{m}(z-y)}.
\end{multline*}
$p_{m}(z-y)\leq C(b)m^{-1},$ $p_{j}(x-y)\geq Cm^{-1}$ for $m/4\leq j\leq m/2.
$ Let $r\overset{\mathrm{def}}{=}m-j,$ which for the region of summation is in
$[m/2,3m/4]$, and $m/2-j-1\leq r/2.$ Then
\begin{align*}
& \mathrm{P}_{0}\left(  X_{1}\neq0, \ldots, X_{m/2-j-1}\neq0, X_{m-j}%
=z-x\right) \\
& \hspace*{2cm} \geq\mathrm{P}_{0}\left(  X_{1}\neq0,\ldots,X_{r/2}%
\neq0,d(X_{r/2},0)\leq\sqrt{m}\right)  \inf_{u:d(u,0)\leq\sqrt{m}}%
\mathrm{P}_{0}(X_{r/2}=z-x-u)\\
& \hspace*{2cm} \geq\frac{C(b)}{m\log m}.
\end{align*}
Therefore, we get
\[
\mathrm{P}_{y,z}^{(m)}\left(  x=X_{j}\,\mathrm{for\,some\,}j\in\left\{
1,\ldots,m\right\}  \right)  \geq\frac{C(b)}{\log m}\geq\frac{C(b)}{\log n}.
\]
\end{proof}

We set
\[
Z_{n,\delta}\overset{\mathrm{def}}{=}\delta\#\left\{  i:\xi_{i}\geq\frac{1}
{4}\log n\right\}  ,
\]
and
\[
\overline{Z}_{n,\delta}\overset{\mathrm{def}}{=}\delta\#\left\{
i:\overline{\xi}_{i}\geq\frac{1}{2}\log n\right\}  .
\]
Then
\[
\mathrm{E}_{\mathbf{X}}\Psi(\mathbf{V}_{b})\geq Z_{n,\delta}\frac{n}{\log
n}\left(  1-\exp\left[  -\frac{c_{4}(b)}{4}\right]  \right)  ,
\]
and therefore
\[
\mathrm{P}_{0}\left(  \mathrm{E}_{\mathbf{X}}\Psi(\mathbf{V}_{b})\leq
\kappa\frac{n}{\log n}\right)  \leq\mathrm{P}_{0}\left(  Z_{n.\delta}
\leq\frac{\kappa}{1-\exp\left[  -c_{4}(b)/4\right]  }\right)  .
\]
Remark now that if $\overline{Z}_{n,\delta}-Z_{n,\delta}\geq8\exp\left[
-c_{1}b^{2}\right]  ,$ then $\Gamma_{n,b,\delta}\geq2\exp\left[  -c_{1}
b^{2}\right]  \frac{\log n}{\delta}.$ Therefore, using (\ref{EstGamma}), we
get
\begin{align*}
\mathrm{P}_{0}\left(  \mathrm{E}_{\mathbf{X}}\Psi(\mathbf{V}_{b})\leq
\kappa\frac{n}{\log n}\right)   &  \leq\mathrm{P}_{0}\left(  \overline
{Z}_{n,\delta}\leq\frac{\kappa}{1-\exp\left[  -c_{4}(b)/4\right]  }
+8\exp\left[  -c_{1}b^{2}\right]  \right) \\
&  +\exp\left[  -c_{2}\frac{\log n}{\delta}\exp\left[  -c_{1}b^{2}\right]
\right]  .
\end{align*}
Choosing $b$ large enough, and then $\kappa>0$ small enough (and
correspondingly $\delta=\kappa^{3}),$ we see that in order to finish the poof
of Proposition \ref{PropModerate}, it suffices to prove the following

\begin{lemma}
For any $R>0$ there exists $\eta>0$ such that for any $\delta>0$
\[
\mathrm{P}_{0}\left(  \overline{Z}_{n,\delta}\leq\eta\right)  \leq n^{-R}%
\]
for $n$ large enough.
\end{lemma}

\begin{proof}
We rescale the random walk by defining
\[
Y_{j}^{(n,\delta)}\overset{\mathrm{def}}{=}X_{jm}/L_{n}.
\]
This random walk depends on $\delta$ through $m=\left[  \delta\frac{n}{\log
n}\right]  .$ It takes values in $T_{n}/L_{n}$ which we regard as a (discrete)
subset of the continuous torus $T\overset{\mathrm{def}}{=}[0,1)^{2}$ with
lattice spacing $1/L_{n}.$ Remember the setting $L_{n}\overset{\mathrm{def}%
}{=}\left[  \sqrt{\frac{n}{\log n}}\right]  $. The transition probabilities of
the $Y$-chain are given by $\tilde{p}(x)=p^{m}(L_{n}x),$ $x\in T_{n}/L_{n}.$
Here $p^{m}$ is the $m$-th matrix power. By the local central limit theorem
(and our aperiodicity assumption) we have that for any $\kappa>0$ there exists
$\gamma_{0}>0$ such that for $\gamma\geq\gamma_{0}$ and any $x\in T_{n}%
/L_{n},$ $n\in\mathbb{N}$ and $\delta>0$
\begin{equation}
\tilde{p}^{\left[  \gamma/\delta\right]  }(x)\leq2L_{n}^{-2}.\label{NearLeb}%
\end{equation}
We denote by $\mathcal{S}_{\delta,\eta}$ the set of unions of square
$[k_{1}\sqrt{\delta},(k_{1}+1)\sqrt{\delta})\times\lbrack k_{2}\sqrt{\delta
},(k_{2}+1)\sqrt{\delta})\subset T$ with total area at most $\eta.$ In order
to prove the lemma, it suffices to prove that for any $R>0$
\begin{equation}
\limsup_{n\rightarrow\infty}\frac{1}{\log n}\log\mathrm{P}_{0}\left(
\sum_{j=0}^{\left(  \log n\right)  /\delta}1_{A}\left(  Y_{j}^{(n,\delta
)}\right)  \geq\frac{\log n}{2\delta}\right)  \leq-R,\label{Superexpon}%
\end{equation}
for small enough $\eta$ uniformly in $\delta$ and $A\in\mathcal{S}%
_{\delta,\eta}.$ We estimate the above probability in a standard way. For any
$\lambda>0$ we have
\begin{equation}
\mathrm{P}_{0}\left(  \sum_{j=0}^{\log n/\delta}1_{A}\left(  Y_{j}%
^{(n,\delta)}\right)  \geq\frac{\log n}{2\delta}\right)  \leq\exp\left[
-\lambda\log n\right]  \mathrm{E}_{0}\left(  \exp\left[  2\lambda\delta
\sum_{j=0}^{\left(  \log n\right)  /\delta}1_{A}\left(  Y_{j}^{(n,\delta
)}\right)  \right]  \right)  .\label{Markov}%
\end{equation}
In order to estimate the right hand side, we use (\ref{NearLeb}). We split the
summation on $j$ alternatively in intervals of length $\gamma/\delta$ and
$3\gamma/\delta,$ the former being called ``short'' intervals, the others
``long''. We begin with a short interval. Remark that the contribution of all
short intervals to the exponent in the expectation on the r.h.s. of
(\ref{Markov}) is at most $\frac{\lambda\log n}{2}.$ Therefore, we can leave
this part out, replacing the first factor on the r.h.s. of (\ref{Markov}) by
$\exp\left[  -\frac{\lambda\log n}{2}\right]  .$ If we choose $\gamma
\overset{\mathrm{def}}{=}\max\left(  \gamma_{0},\frac{\log2}{\lambda}\right)
$ we have by (\ref{NearLeb})
\begin{multline*}
\mathrm{E}_{0}\left(  \exp\left[  2\lambda\delta\sum_{j\in
\mathrm{long\,\,intervals}}1_{A}\left(  Y_{j}^{(n,\delta)}\right)  \right]
\right) \\ \leq\exp\left[  \frac{\lambda\log n}{4}\right]  \left\{
\mathrm{E}_{u}\left(  \exp\left[  2\lambda\delta\sum_{j=0}^{3\gamma/\delta
}1_{A}\left(  Y_{j}^{(n,\delta)}\right)  \right]  \right)  \right\}
^{\frac{\log n}{4\gamma}},
\end{multline*}
where $\mathrm{E}_{u}$ is the expectation with respect to an uniform starting
distribution. We therefore get
\begin{align*}
&  \limsup_{n\rightarrow\infty}\frac{1}{\log n}\log\mathrm{P}_{u}\left(
\sum_{j=0}^{\left(  \log n\right)  /\delta}1_{A}\left(  Y_{j}^{(n,\delta
)}\right)  \geq\frac{\log n}{2\delta}\right) \\
&  \leq-\frac{\lambda}{4}+\frac{1}{4\gamma}\lim_{n\rightarrow\infty}%
\log\mathrm{E}_{u}\left(  \exp\left[  2\lambda\delta\sum_{j=0}^{3\gamma
/\delta}1_{A}\left(  Y_{j}^{(n,\delta)}\right)  \right]  \right) \\
&  =-\frac{\lambda}{4}+\frac{1}{4\gamma}\log\mathrm{E}_{u}\left(  \exp\left[
2\lambda\delta\sum_{j=0}^{3\gamma/\delta}1_{A}\left(  B_{\delta j}\right)
\right]  \right)  ,
\end{align*}
where $\left(  B_{t}\right)  _{t\geq0}$ is a Brownian motion on $T$ with
covariance matrix $\mathcal{Q}.$ For $x\geq0$ we have \textrm{e}$^{x}%
\leq1+x\mathrm{e}^{x},$ and we therefore get
\begin{align*}
\mathrm{E}_{u}\left(  \exp\left[  2\lambda\delta\sum_{j=0}^{3\gamma/\delta
}1_{A}\left(  B_{\delta j}\right)  \right]  \right)   &  \leq1+2\lambda
\delta\sum_{j=0}^{3\gamma/\delta}\mathrm{P}_{u}\left(  B_{\delta j}\in
A\right)  \mathrm{e}^{6\lambda\gamma}\\
&  =1+6\lambda\gamma\left|  A\right|  \mathrm{e}^{6\lambda\gamma}%
\leq1+6\lambda\gamma\eta\mathrm{e}^{6\lambda\gamma}.
\end{align*}
We therefore get
\[
\limsup_{n\rightarrow\infty}\frac{1}{\log n}\log\mathrm{P}_{u}\left(
\sum_{j=0}^{\left(  \log n\right)  /\delta}1_{A}\left(  Y_{j}^{(n,\delta
)}\right)  \geq\frac{\log n}{2\delta}\right)  \leq-\frac{\lambda}{4}%
+\frac{3\lambda\eta}{2}\mathrm{e}^{6\lambda\gamma}.
\]
Choosing $\lambda$ appropriately, this proves the claim
\end{proof}

\section{The case d=1}

\label{app_1d}

We consider the $\delta$-pinning case only, and $p\left(  \pm1\right)  =1/2.$
We however can easily allow more general symmetric interaction functions
$V:\mathbb{R}\rightarrow\mathbb{R}^{+}.$ We set
\[
\psi\left(  x\right)  \overset{\mathrm{def}}{=}\frac{\mathrm{e}^{-\beta
V(x)/2}}{\int\mathrm{e}^{-\beta V(y)/2}\,\mathrm{d}y}.
\]
The only property we need is $\int\mathrm{e}^{-\beta V(y)/2}\,\mathrm{d}%
y<\infty,$ $\int x\psi\left(  x\right)  \,\mathrm{d}x=0,$ $\int x^{2}%
\psi\left(  x\right)  \,\mathrm{d}x=\sigma^{2}<\infty.$ By a simple rescaling,
we can assume $\sigma^{2}=1.$ Let $\psi_{k}$ be the $k$-fold convolution of
$\psi.$ By the local central limit theorem, we have
\[
f\left(  k\right)  \overset{\mathrm{def}}{=}\psi_{k}\left(  0\right)
=\frac{1}{\sqrt{2\pi k}}+o\left(  \frac{1}{\sqrt{k}}\right)  ,
\]
as $k\rightarrow\infty.$

The distribution $\nu_{n}^{\varepsilon}$ of pinned sites on $\Lambda=\left\{
-n,-n+1,\ldots,n\right\}  $ is easily described: Let $A\subset\left\{
-n,-n+1,\ldots,n\right\}  $ with $\left|  A\right|  =m-1,$%
\[
A=\left\{  k_{1},k_{2},\ldots,k_{m-1}\right\}  ,
\]
where $k_{0}\overset{\mathrm{def}}{=}-n-1<k_{1}<k_{2}<\ldots<k_{m-1}%
<k_{m}\overset{\mathrm{def}}{=}n+1.$ Then
\begin{equation}
\nu_{n}^{\varepsilon}\left(  A\right)  =\frac{1}{Z_{n,\varepsilon}}%
\varepsilon^{m-1}\prod_{j=1}^{m}f\left(  k_{j}-k_{j-1}\right)
.\label{renewal}%
\end{equation}
Of course, $\sum_{k}f(k)=\infty.$ Therefore, there exists a unique
$\lambda=\lambda\left(  \varepsilon\right)  ,$ such that
\[
\varepsilon\sum_{k}\mathrm{e}^{-\lambda k}f(k)=1.
\]
Remark that (\ref{renewal}) is not changed if we replace $f$ by
\[
f_{\lambda}(k)\overset{\mathrm{def}}{=}\mathrm{e}^{-\lambda k}f(k).
\]
Standard renewal arguments then show that $\nu^{\varepsilon}=\lim
_{n\rightarrow\infty}\nu_{n}^{\varepsilon}$ exists, and is simply given as the
stationary renewal sequence with renewal epochs with distribution $\left\{
f_{\lambda\left(  \varepsilon\right)  }(k):k>0\right\}  .$ For instance, if
\begin{align*}
& \xi\overset{\mathrm{def}}{=}\max\left\{  m\leq0:m\in\mathcal{A}\right\} \\
& \eta\overset{\mathrm{def}}{=}\min\left\{  m>0:m\in\mathcal{A}\right\}  ,
\end{align*}
then

\begin{lemma}%
\[
\nu^{\varepsilon}\left(  \left(  \xi,\eta\right)  =\left(  k,l\right)
\right)  =\frac{1}{M^{\varepsilon}}f_{\lambda\left(  \varepsilon\right)
}(l-k)
\]
if $k\leq0<l,$ where $M^{\varepsilon}\overset{\mathrm{def}}{=}\sum
_{j}j\,f_{\lambda\left(  \varepsilon\right)  }\left(  j\right)  .$
\end{lemma}

The full measure $\mu^{\varepsilon}$ (in the thermodynamic limit) is then
given as a mixture
\[
\mu^{\varepsilon}=\sum_{A\subset\mathbb{Z}}\nu^{\varepsilon}\left(  A\right)
\mu_{A},
\]
where $\mu_{A}$ is the measure on $\mathbb{R}^{\mathbb{Z}}$ given by
independent pieces of tied-down random walks between successive elements of
$A.$ For instance
\[
\int\phi_{0}^{2}\,\mu^{\varepsilon}\left(  d\phi\right)  =\sum_{k\leq0<l}%
\frac{1}{M^{\varepsilon}}f_{\lambda\left(  \varepsilon\right)  }%
(l-k)\mathrm{E}_0\left(  S_{-k}^{2}|S_{l-k}=0\right)  ,
\]
where $S_{0},S_{1},S_{2},\ldots$ is a random walk on $\mathbb{R}$ starting at
$0$ with
distribution of the increments given by $\psi.$ We now want to determine the
$\varepsilon\rightarrow0$ behavior of this quantity. First, remark that for
small $\lambda>0$%
\begin{align*}
\varepsilon\sum_{k}\mathrm{e}^{-\lambda k}f(k)  & \sim\varepsilon\sum
_{k}\mathrm{e}^{-\lambda k}\frac{1}{\sqrt{2\pi k}}\sim\frac{\varepsilon}%
{\sqrt{2\pi\lambda}}\sum_{k}\lambda\frac{1}{\sqrt{k\lambda}}\mathrm{e}%
^{-\lambda k}\\
& \sim\frac{\varepsilon}{\sqrt{2\pi\lambda}}\int_{0}^{\infty}\frac{1}{\sqrt
{x}}\mathrm{e}^{-x}\,\mathrm{d}x=\frac{\varepsilon}{\sqrt{2\lambda}}.
\end{align*}
Therefore
\[
\lambda\left(  \varepsilon\right)  =\frac{\varepsilon^{2}}{2}+o\left(
\varepsilon^{2}\right)  .
\]
From this, we get%
\begin{align}
M_{\varepsilon}  & =\sum_{j}j\,f_{\lambda\left(  \varepsilon\right)  }\left(
j\right)  \sim\sum_{j}j\,\frac{1}{\sqrt{2\pi j}}\mathrm{e}^{-\varepsilon
^{2}j/2}\nonumber\\
& =\frac{1}{\varepsilon^{3}\sqrt{2\pi}}\sum_{j}\varepsilon^{2}\sqrt
{\varepsilon^{2}j}\mathrm{e}^{-\varepsilon^{2}j/2}\label{Mean}\\
& \sim\frac{1}{\varepsilon^{3}\sqrt{2\pi}}\int_{0}^{\infty}\sqrt{x}%
\mathrm{e}^{-x/2}\,\mathrm{d}x=\frac{1}{\varepsilon^{3}}.\nonumber
\end{align}
Furthermore%
\begin{align*}
\sum_{k\leq0<l}f_{\lambda\left(  \varepsilon\right)  }(l-k)\mathrm{E}\left(
S_{-k}^{2}|S_{l-k}=0\right)    & \sim\sum_{k\leq0<l}\frac{1}{\sqrt{2\pi\left(
l-k\right)  }}\mathrm{e}^{-\varepsilon^{2}\left(  l-k\right)  /2}%
\mathrm{E}_{l-k}\left(  S_{-k}^{2}\right)  \\
& =\sum_{n=1}^{\infty}\frac{1}{\sqrt{2\pi n}}\mathrm{e}^{-\varepsilon^{2}%
n/2}\sum_{m=0}^{n-1}\mathrm{E}_{n}\left(  S_{m}^{2}\right)
\end{align*}
where $\mathrm{E}_{m}$ stands for the expectation with respect to a random
walk tied down after time $m.$ The right-hand side of the above expression is%
\begin{align*}
& \sim\sum_{n=1}^{\infty}\frac{\sqrt{n}}{\sqrt{2\pi}}\mathrm{e}^{-\varepsilon
^{2}n/2}\sum_{m=0}^{n-1}\frac{m}{n}\left(  1-\frac{m}{n}\right)  \\
& \sim\frac{1}{6\varepsilon^{5}}\frac{1}{\sqrt{2\pi}}\sum_{n=1}^{\infty
}\varepsilon^{2}\left(  \varepsilon^{2}n\right)  ^{3/2}\mathrm{e}%
^{-\varepsilon^{2}n/2}\sim\frac{1}{6\varepsilon^{5}}\frac{1}{\sqrt{2\pi}}%
\int_{0}^{\infty}y^{3/2}\,\mathrm{e}^{-y/2}\,\mathrm{d}y=\frac
{1}{2\varepsilon^{5}}.
\end{align*}
Combining this with (\ref{Mean}) yields

\begin{proposition}%
\[
\int\phi_{0}^{2}\,\mu^{\varepsilon}(\mathrm{d}\phi)=\frac{1}{2\varepsilon^{2}%
}+o\left(  \frac{1}{\varepsilon^{2}}\right)
\]
as $\varepsilon\rightarrow0.$
\end{proposition}

The mass is very easy, too. For fixed $\varepsilon$, the $x\rightarrow\infty$
limit of $\int\phi_{0}\phi_{x}\,\mu^{\varepsilon}(\mathrm{d}\phi)$ is in
leading order the same as the probability under $\nu^{\varepsilon}$ that the
interval $\left[  0,x\right]  $ has no renewal point. In leading order, this
is just the exponential tail behavior of the distribution $f_{\lambda\left(
\varepsilon\right)  }.$ Therefore, we get

\begin{proposition}%
\[
\lim_{x\rightarrow\infty}\frac{1}{x}\log\int\phi_{0}\phi_{x}\,\mu
^{\varepsilon}(\mathrm{d}\phi)=-\lambda\left(  \varepsilon\right)
=-\frac{\varepsilon^{2}}{2}+o\left(  \varepsilon^{2}\right)  .
\]
\end{proposition}

\vfill

\noindent\textsc{E. Bolthausen, Institut f\"ur Mathematik, Universit\"{a}t
Z\"{u}rich, Winterthurerstrasse 190, CH--8057 Z\"urich, Switzerland}%
\newline E--mail address: \texttt{eb@amath.unizh.ch}\newline Homepage:
\texttt{http://www.math.unizh.ch/eb/}

\vskip0.3 cm \noindent\textsc{Y. Velenik, Faculty of industrial Engineering,
Technion, Haifa 32000, Israel\newline New address: L.A.T.P., UMR-CNRS 6632,
C.M.I., 39 rue Joliot Curie, 13453 Marseille cedex 13, France.}\newline
E--mail address: \texttt{velenik@cmi.univ-mrs.fr}\newline Homepage:
\texttt{http://www.cmi.univ-mrs.fr/$\scriptstyle\sim$velenik/}
\end{document}